\documentclass[11pt,a4paper,leqno]{article}

\usepackage{amsmath,amscd,amsfonts,amsthm,amssymb,dsfont,bbm,manfnt}
\usepackage{url,color,appendix,eucal,comment,tabu,float}
\usepackage[usenames,dvipsnames]{xcolor}
\usepackage[nottoc,numbib]{tocbibind}
\usepackage{multicol}
\usepackage[OT2,T1]{fontenc}
\usepackage[utf8]{inputenc}
\usepackage[all]{xy}
\usepackage{lscape,rotating} 
\usepackage{graphicx}
\usepackage[nolabel]{showlabels}
\usepackage[hypertexnames=false]{hyperref} 

\newenvironment{pf}
{\medskip\noindent {\it Proof.  }}
{\hfill\nobreak $\Box$ \par\bigbreak}

\newcommand{\ps}{\par \smallskip}

\newcommand{\Z}{\mathbb{Z}}

\newcommand{\R}{\mathbb{R}}

\newcounter{introcounter}
\newtheorem{thm}[subsection]{Theorem}
\newtheorem{lemma}[subsection]{Lemma}

\newtheorem{remark}[subsection]{Remark}

\newtheorem{prop}[subsection]{Proposition}
\newtheorem{example}[subsection]{Example}

\newtheorem{thmintro}{Theorem}

\newtheorem{definition}[subsection]{Definition}

\newtheorem{fact}[subsection]{Fact}

\usepackage{titlesec}

\titleformat{\section}{\bf \normalsize}{\arabic{section}.}{1 em}{}

\titleformat{\subsection}
{}
{\bf \arabic{section}.\arabic{subsection}.}
{0.5 em}
{\normalsize}

\titleformat{\subsubsection}[runin]
{\small \bf}
{}
{}
{}

\hypersetup{colorlinks=true,linkcolor=black,citecolor=black,urlcolor=Periwinkle,pdfborderstyle={/S/U/W 1}}

\makeatletter

\@addtoreset{equation}{section}
\makeatother

\makeatletter

\@addtoreset{table}{section}
\makeatother

\begin{document}

\title{Unimodular hunting II}

\author{
Bill Allombert\thanks{\texttt{Bill.Allombert@math.u-bordeaux.fr}, CNRS, INRIA, Institut de Math\'ematiques de Bordeaux, Universit\'e de Bordeaux, 351 cours de la Lib\'eration, F 33 405 Talence, France.}\\ 
\and 
Ga\"etan Chenevier\thanks{\texttt{gaetan.chenevier@math.cnrs.fr}, CNRS, \'Ecole Normale Sup\'erieure-PSL, D\'epartement de Math\'ematiques et Application, 45 rue d'Ulm, 75230 Paris Cedex, France. During this work, G. Chenevier was supported by the project ANR-19-CE40-0015-02 COLOSS.}
}


\maketitle

\begin{abstract}
Pursuing ideas in \cite{cheuni}, we determine the isometry classes of unimodular lattices of rank $28$,  
as well as the isometry classes of unimodular lattices of rank $29$ without nonzero vectors of norm $\leq 2$.
\end{abstract}

\ps
\ps

\section{Introduction}
\label{sect:intro}
${}^{}$ 
For $n\geq 1$, we denote by $\mathcal{L}_n$ the set of unimodular integral lattices in the standard Euclidean space $\R^n$, and by ${\rm X}_n$ the set of isometry classes of elements of $\mathcal{L}_n$. In the recent work \cite{cheuni}, the second author studied the cyclic Kneser neighbors of the standard lattice ${\rm I}_n:=\Z^n$ and classified the elements of ${\rm X}_n$ for $n=26$ and $27$, developing a method initiated by Bacher and Venkov \cite{bacher, bachervenkov}.
Let ${\rm X}_n^\emptyset \subset {\rm X}_n$ be the subset of classes of lattices $L \in \mathcal{L}_n$ with no nonzero vector $v \in L$ such that $v.v \leq 2$. The main result of this paper is the following.

\begin{thmintro} \label{X28}
We have $|{\rm X}_{28}|=374062$ and $|{\rm X}^\emptyset_{29}|=10092$. 
\end{thmintro}

We refer to~\cite{cheuni} for some historical background on these questions, especially in lower dimensions.
As we have $|{\rm X}_{27}|=17059$, the assertion about $|{\rm X}_{28}|$ is equivalent to $|{\rm X}'_{28}|=357003$, where ${\rm X}'_n \subset {\rm X}_n$ denotes the subset consisting of classes of lattices with no norm $1$ vectors.
The $38$ lattices in ${\rm X}_{28}^\emptyset$ had already been determined by Bacher and Venkov in \cite{bachervenkov}. Moreover, King's refinements of the Minkowski-Siegel-Smith mass formulae in \cite{king} gave the (not too far!) lower bounds $|{\rm X}_{28}|> 327972$ and $|{\rm X}_{29}^\emptyset|> 8911$. \ps 

As in \cite{cheuni}, our aim is not only to determine ${\rm X}_{28}$ and ${\rm X}^\emptyset_{29}$ but also to provide representatives of each isometry class as a (cyclic) $d$-neighbor of the standard lattice ${\rm I}_n:=\Z^n$. We briefly recall their concrete definition and refer to {\it loc. cit.} for more details. Let $d\geq 1$ and $x \in \Z^n$ with ${\rm gcd}(x_1,\dots,x_n)=1$. Then ${\rm M}_d(x) :=\{ v \in \Z^n\,\,|\,\,\sum_{i=1}^n v_i x_i \equiv 0 \bmod d\}$ is an index $d$ sublattice of ${\rm I}_n$. Setting $e=1$ for $d$ odd, $e=2$ otherwise, and assuming furthermore $x$ is $d$-{\it isotropic}, that is $\sum_{i=1}^n x_i^2 \equiv 0 \bmod ed$, then there are exactly $e$ unimodular lattices $L \in \mathcal{L}_n$
with $L \cap {\rm I}_n = {\rm M}_d(x)$; they have the form
$${\rm N}_d(x')\,:= \,{\rm M}_d(x) \,+\, \frac{x'}{d} \Z$$
where $x' \in \Z^n$ satisfies $x'.x' \equiv 0 \bmod d^2$ and $x' \equiv x \bmod d$. 
For $d$ odd, this unique lattice is denoted ${\rm N}_d(x)$ or ${\rm N}_d(x; 0)$. For $d$ even, 
it only depends on the $\epsilon \in \{ 0,1\}$ defined by $2 x'.x \equiv x.x + \epsilon d^2 \bmod 2d^2$, 
and is denoted ${\rm N}_d(x; \epsilon)$. 

\begin{thmintro} 
\label{listdxe}
A list of $(d,x,\epsilon)$ such that the ${\rm N}_d(x; \epsilon)$ are representatives for the isometry classes in ${\rm X}'_{28}$ and ${\rm X}_{29}^\emptyset$ are given \href{http://gaetan.chenevier.perso.math.cnrs.fr/unimodular_lattices/L28}{in}~{\rm \cite{allcheweb28}} and~\href{http://gaetan.chenevier.perso.math.cnrs.fr/unimodular_lattices/unimodular_lattices_no_roots.gp}{in}~{\rm \cite{allcheweb29}}.
\end{thmintro}

	Let us first comment on these lists, starting with the $10092$ elements in ${\rm X}_{29}^\emptyset$.
The statistics for the order of their isometry groups are given in Table~\ref{tab:roaX29empty}. In particular, 
about $80$ \% of them have trivial isometry group $\{ \pm 1\}$.

\begin{table}[H]
\tabcolsep=3pt
{\scriptsize \renewcommand{\arraystretch}{1.7} \medskip
\begin{center}
\begin{tabular}{c|c|c|c|c|c|c|c|c|c|c|c|c|c}
\texttt{ord} & $2$ &  $4$  & $6$ & $ 8$ & $12$ & $16$ & $20$ & $24$ & $32$ & $36$ & $40$ &  $48$ & $60$  \\
$\sharp$ & $8081$ & $1465$ & $6$ & $293$ & $28$ & $91$ & $1$ & $21$ & $32$ & $1$ & $3$ & $15$ & $1$  \\
\hline
\texttt{ord} & $64$ & $72$ & $80$ & $96$ & $120$ & $128$ & $144$ & $160$ & $192$ & $232$ & $256$ & $288$ & $320$    \\ 
$\sharp$ & $12$ & $1$ & $1$ & $11$ & $1$ & $2$ & $1$ & $2$ & $2$ & $1$ & $2$ & $1$ & $1$ \\
 \hline
 \texttt{ord} & $384$ & $768$ & $864$ & $960$ & $1024$ & $1536$ & $2400$ & $2592$ & $3072$ & $5184$ & $6144$ & $18432$ & $24000$  \\
 $\sharp$ & $1$ & $2$ &  $1$ & $1$ & $2$ & $2$ & $1$ & $1$ & $1$ & $1$ & $1$ & $1$ & $1$
\end{tabular} 
\end{center}
} 
\caption{Number $\sharp$ of classes of lattices $L$ in ${\rm X}_{29}^\emptyset$ with $|{\rm O}(L)|=\texttt{ord}$.}
\label{tab:roaX29empty} 
\end{table}

Arguing as in~\cite[Sect. 12]{cheuni} (in particular, using \cite{GAP}), we checked that
only two of these $10092$ lattices have a non-solvable isometry group\footnote{Actually, there are only five $L$ 
such that $|{\rm O}(L)|$ is both $\equiv 0 \bmod 4$ and not of the form $p^a q^b$ with $p,q$ primes.}. This is very little compared to case of lower dimensions (see {\it loc. cit.}).
Their isometry groups have order $2400$ and $960$, and are 
respectively isomorphic to $\Z/2 \times {\rm S}_5 \times {\rm D}_{10}$ and $\Z/2 \times \Z/4 \times {\rm S}_5$. These two lattices are furthermore {\it exceptional} in the sense of Bacher and Venkov, {\it i.e.} have a characteristic vector of norm $5$
(see \S \ref{sect:notation} for the unexplained terminology in this introduction). Actually, exactly $105$ lattices in ${\rm X}_{29}$  are exceptional, and the statistics for the order of their isometry groups are given by Table~\ref{tab:roaX29emptyex}.
\begin{table}[H]
\tabcolsep=3pt
{\scriptsize \renewcommand{\arraystretch}{1.8} \medskip
\begin{center}
\begin{tabular}{c|c|c|c|c|c|c|c|c|c|c|c|c|c|c|c|c|c}
\texttt{ord} & $2$ & $4$ & $8$ & $12$ & $16$ & $20$ & $32$ & $40$ & $48$ & $64$ & $80$ & $96$ & $192$ & $960$ & $1024$ & $2400$ & $24000$  \\
$\sharp$ & $20$ & $31$ & $24$ & $2$ & $10$ & $1$ & $3$ & $3$ & $3$ & $1$ & $1$ & $1$ & $1$ & $1$ & $1$ & $1$ & $1$ 
\end{tabular} 
\end{center}
} 
\caption{Number $\sharp$ of exceptional lattices $L$ in ${\rm X}_{29}^\emptyset$ with $\texttt{ord}=|{\rm O}(L)|$.}
\label{tab:roaX29emptyex} 
\end{table}

	We now consider the $357003$ rank $28$ unimodular lattices with no norm $1$ vectors. 
According to King, there are $4722$ possible root systems for them: see \cite{cheuni} Table 1.2. 
We found it convenient\footnote{Almost all unimodular lattices $L$ in our lists with ${\rm i}(L)=i$
are of the form ${\rm N}_d(x;\epsilon)$ where the pair $(d,x)$ has index $i$ in the sense of Sect. \ref{sect:thebiasedalgo}.}
to split these lattices $L$ according to the integer ${\rm i}(L)$ defined as the maximal integer $i\geq 1$ such that the root lattice ${\rm A}_{i-1}$ may be embeded into $L$. For instance, we have ${\rm i}(L)=1$ (resp. $2$, $3$) if, and only if, the root system of $L$ is empty (resp. $r {\bf A}_1$, resp. $r{\bf A}_2\, \, s {\bf A}_1$ with $r\geq 1$). 

\begin{table}[H]
\tabcolsep=4pt
{\scriptsize \renewcommand{\arraystretch}{1.8} \medskip
\begin{center}
\begin{tabular}{c|c|c|c|c|c|c|c|c|c|c|c|c|c|c}
\texttt{i} & $1$ & $2$ & $3$ & $4$ & $5$ & $6$ & $7$ & $8$ & $9$ & $10$ & $11$ & $12$ & $13$ & $14$ \\
$\sharp$ & $38$ & $20560$ & $121684$ & $126661$ & $55585$ & $20919$ & $6712$ & $2935$ & $960$ & $516$ & $168$ & $142$ & $45$ & $35$ \\  \hline
\texttt{i} & $15$ & $16$ & $17$ & $18$ & $19$ & $20$ & $21$ & $22$ & $23$ & $24$ & $25$ & $26$ & $27$ & $28$ \\
$\sharp$ &  $8$ & $20$ & $3$ & $3$ & $1$ & $5$ & $1$ & $0$ & $0$ & $1$ & $0$ & $0$ & $0$ & $1$
\end{tabular} 
\end{center}
} 
\caption{Number $\sharp$ of lattices $L$ in ${\rm X}'_{28}$ with $\texttt{i}={\rm i}(L)$.}
\label{tab:numbX28giveni} 
\end{table}
Only $238$ of the $357003$ lattices in ${\rm X}'_{28}$ have a non solvable {\it reduced isometry group} (see \S 
\ref{sect:notation} below). Four of them have a Jordan-H\"older factor not appearing for unimodular lattices of smaller rank. Those $4$ lattices have an empty root system, and thus belong the Bacher-Venkov list; their isometry groups are described in Table~\ref{tab:X28newsimpl}.

\tabcolsep=3.5pt
\begin{table}[H]
{\scriptsize 
\renewcommand{\arraystretch}{1.8} \medskip
\begin{center}
\begin{tabular}{c|c|c|c|c}
$|G|$ & $9170703360$ & $348364800$ & $4838400$ & $58240$  \\  \hline
$G$ &  \scalebox{.9}{$(\Z/2 \,. \,{\rm PSp}_6(3)) : \Z/2$} &  \scalebox{.9}{$\Z/2 \times ({\rm O}^+_8(2) : \Z/2)$} &  \scalebox{.9}{$(((\Z/2 \,.\, {\rm HJ}) : \Z/2) : \Z/2) : \Z/2$} &  \scalebox{.9}{$\Z/4 \times {\rm Sz}(8)$}
\end{tabular} 
\end{center}
}
\caption{{\small The isometry groups $G$ of the $4$ lattices in ${\rm X}'_{28}$ whose reduced isometry group has a ``new'' Jordan-H\"older factor, using \texttt{GAP}'s notations.}}
\label{tab:X28newsimpl}
\end{table}

	We finally discuss the proofs of Theorems~\ref{X28} and~\ref{listdxe}. 
The ingredients are the same as those described in the introduction of~\cite{cheuni}: 
systematic study of all the cyclic Kneser $d$-neighbors of ${\rm I}_n$ for $d=2, 3, ...$ ("coupon-collector'' problem), bet on fine enough isometry invariants, 
splitting according to root systems thanks to King's results \cite{king}, use of clever visible root systems to biase the search, and case-by-case more specific methods for the remaining lattices with small masses (suitable $2$-neighbors, exceptional lattices, ``addition of ${\rm D}_m$'', visible isometries). \par

However, compared to the work {\it loc. cit.}, 
the computations here are of a much larger scale, as the number of lattices in Theorem~\ref{X28} already indicates. 
There are in particular $4\,722$ different root systems $R$ such that ${\rm X}_{28}^R$ is non empty, and during our search we were forced to study several hundreds of them case by case. In the last Sect.~\ref{sect:examples} we give five examples  of such a study (including $R=\emptyset$ in dimension $29$) in order to illustrate some important techniques that we 
used. Let us insist now on a few novelties that made our computations possible.\ps\ps

{\bf 1. The choice of invariant}. It is crucial for our method to have at our disposal 
an invariant which is both fine enough to distinguish all of our lattices, 
and fast to compute. For a lattice $L$ and $i\geq 0$, 
we set 
{\small 
\begin{equation}\label{eq:defri}
{\rm R}_{\leq i}(L)=\{ v \in L\, |\, v.v \leq i\},\,\, {\rm R}_{i}(L)=\{ v \in L\, |\, v.v = i\}\,\,{\rm 
and} \,\,{\rm r}_i(L)=|{\rm R}_{i}(L)|.
\end{equation}
}
Unimodular lattices of rank $\leq 29$ tend to be generated over $\Z$ by their ${\rm R}_{\leq 3}$, 
which is thus a natural candidate for an invariant. 
However, we are not aware of any classification of finite metric sets 
of the form ${\rm R}_{\leq 3}(L)$ or ${\rm R}_3(L)$, 
contrary to the case of ${\rm R}_{\leq 2}(L)$ (which is the theory of root systems). 
In \S\ref{sect:inv} we define, for an integral lattice $L$, 
an invariant that we denote ${\rm BV}(L)$ and 
which only depends on ${\rm R}_{\leq 3}(L)$. It is a variant of one of the invariants 
used by Bacher and Venkov in their classification of 
${\rm X}_{28}^\emptyset$ and ${\rm X}_{27}^\emptyset$, 
hence the notation. A surprising fact, which eventually follows from our computation, 
but for which we do not yet have a theoretical explanation, is the following.

\begin{thmintro} 
\label{thmi:invnorm3}
Let $L$ and $L'$ be two rank $n$ unimodular lattices with ${\rm r}_1(L)={\rm r}_1(L')=0$.  Assume either $n\leq 28$, or $n=29$ and ${\rm r}_2(L)={\rm r}_2(L')=0$. Then $L$ and $L'$ are isometric if, and only if, we have ${\rm BV}(L)={\rm BV}(L')$.
\end{thmintro}

In particular, $L$ and $L'$ are isometric if, and only if, ${\rm R}_{\leq 3}(L)$ and ${\rm R}_{\leq 3}(L')$ are isometric. \ps\ps

{\bf 2. Algorithmic improvements}. 
In \S\ref{sect:Zbasisshort}, we discuss a simple probabilistic algorithm whose aim is to find $\Z$-basis consisting of small vectors of a given lattice, and which allowed to substantially shorten the computation of the order of the (reduced) isometry groups of our lattices. In \S\ref{sect:thebiasedalgo}, we give some details about the Biased Neighbor Enumeration algorithm 
informally described in \cite[\S 1.10]{cheuni} and that we repeatedly used in our search. \par
The programs used to perform the computations were written in the GP language with some critical parts in C using the libpari library. They were run on the two clusters PlaFRIM and Cinaps.
We used the parallel programming interface of GP to make efficient use of the clusters,
by allowing to switch between POSIX threads and MPI depending on the hardware available.\par
The total CPU time was about $72$ years. 
Fortunately, it may be checked independently, and a posteriori, 
that the given lists are complete: it is enough to
check that our lattices have distinct ${\rm BV}$ invariants and that the sum of
their masses coincides with the mass of ${\rm X}_n$. See \cite{allcheweb28} and \cite{allcheweb29} for the relevant \texttt{PARI/GP}
source code for this check. This is much shorter, and ``only'' requires
a few hours in dimension $29$, and about $27$ days in dimension $28$. 
It simultaneously proves Theorems~\ref{X28},~\ref{listdxe} and~\ref{thmi:invnorm3}.
\footnote{As a consequence,
this also provides an independent verification of King’s computations in \cite{king}.}\ps
\ps \medskip
{\sc Acknowledgments.}  The authors thank Jacques Martinet and Olivier Ta\"ibi for remarks or discussions. 
Experiments presented in this work were carried out using: (a) the PlaFRIM experimental testbed, supported by Inria, CNRS (LABRI and IMB),
Universite de Bordeaux, Bordeaux INP and Conseil Regional d'Aquitaine (see 
\url{https://www.plafrim.fr/}), (b)  the cluster cinaps of the LMO, Universit\'e Paris-Saclay. 
We warmly thank PlaFRIM and the LMO for sharing their machines.


\section{Notation and terminology}  
\label{sect:notation}
We use mostly classical notation and terminology. 
Let $L$ be an integral Euclidean lattice. The (finite) isometry group of $L$ is denoted by ${\rm O}(L)$. 
The notations ${\rm R}_{\leq i}(L)$, ${\rm R}_i(L)$ and ${\rm r}_i(L)$ have been introduced in Formula~\eqref{eq:defri}.
We refer {\it e.g.} to~\cite{cheuni} \S 2, \S 4.3 and \S 9 for complements. \ps

-- The {\it norm} of a vector $v \in L$ is $v.v$. The lattice $L$ is even if we have $v.v \in 2\Z$ for all $v \in L$, odd otherwise. A {\it characteristic vector} of $L$ is a vector $\xi \in L$ such that $\xi.v \equiv v.v \bmod 2$ for all $v \in L$. 
For $L \in \mathcal{L}_n$ odd, the norm of such a vector is $\equiv n \bmod 8$, and we denote by
${\rm Exc}(L)$ the set of characteristic vectors of $L$ of norm $<8$.
Following \cite{bachervenkov}\, \S 3, a lattice $L \in \mathcal{L}_n$ is called {\it exceptional} 
if we have ${\rm Exc}(L) \neq \emptyset$. \ps
-- The {\it root system} of $L$ is the finite metric set $R:={\rm R}_2(L)$; this is an ${\rm ADE}$ root system.\footnote{In this papers, 
roots are always assumed to have norm $2$. } The {\it Weyl group} of $L$ is the subgroup ${\rm W}(L) \subset {\rm O}(L)$ generated by the orthogonal symetries about each $\alpha \in {\rm R}_2(L)$. This is a normal subgroup isomorphic to the classical Weyl group ${\rm W}(R)$ of $R$. 
We define the {\it reduced isometry group} of $L$ to be ${\rm O}(L)/{\rm W}(L)$.\ps

-- For an arbitrary ${\rm ADE}$ root system $S$, we denote by ${\rm Q}(S)$ the even lattice it generates ({\it root lattices}).
We use bold fonts ${\bf A}_n$, ${\bf D}_n$, ${\bf E}_n$ to denote 
isomorphism classes of root systems of that names, with the conventions ${\bf A}_0={\bf D}_0={\bf D}_1=\emptyset$ and ${\bf D}_2= 2{\bf A}_1$.  
We also denote by ${\rm A}_n$ ($n\geq 1$), ${\rm D}_n$ ($n\geq 2$) and ${\rm E}_n$ ($n=6,7,8$) the standard corresponding root lattices.\ps

-- The {\it mass} of any collection $\mathcal{L}$ of lattices, denoted ${\rm mass} \,\mathcal{L}$, is defined as the sum, over representatives $L_i$ of the isometry classes of lattices in $\mathcal{L}$, of $1/|{\rm O}(L_i)|$. In particular, if $L$ is a (single) integral euclidean lattice, the mass of $L$ is $1/|{\rm O}(L)|$. We also define the {\it reduced mass} of $L$ as $\frac{|{\rm W}(L)|}{|{\rm O}(L)|}$. \ps
-- We denote by ${\rm X}_n^R \subset {\rm X}_n$ the subset of isometry classes of lattices $L$ 
in $\mathcal{L}_n$ with ${\rm r}_1(L)=0$ and ${\rm R}_2(L) \simeq R$. The {\it reduced mass} of ${\rm X}_n^R$ is also defined as ${\rm rmass} \,{\rm X}_n^R \,:=\,|{\rm W}(R)|\, {\rm mass}\,{\rm X}_n^R$. We have the interesting lower bound $|{\rm X}_n^R| \geq  m\, {\rm rmass}\,{\rm X}_n^R$, with $m=1$ if $R$ has rank $n$ and ${\rm W}(R)$ contains $-1$, and $m=2$ otherwise. 
It follows from \cite{king} that, for all $n\leq 30$ and all $R$, we know the rational ${\rm rmass}\,\,{\rm X}_n^R$ (see also \cite[\S 6]{cheuni}). \ps
-- We finally recall a few specific features of unimodular lattices 
in dimension $n \equiv 4 \bmod 8$ (such as $n=28$), as well as the description of exceptional lattices 
in this case given in \cite[\S 9]{cheuni}. Assume $n \equiv 4 \bmod 8$ and $L\in \mathcal{L}_n$. As is well-known,  
there are exactly two other $L'$ in $\mathcal{L}_n$ having the same {\it even part};\footnote{The even part of an odd integral lattice $L$ is the index $2$ sublattice $\{x \in L\, \, |\, x.x \equiv 0 \bmod 2\}$.}
we call these two $L'$ the {\it companions} of $L$. 
Assume ${\rm r}_1(L)=0$. Then $L$ is exceptional if, and only if, it has a companion $L'$ with ${\rm r}_1(L') \neq 0$.
In this case, this $L'$ is unique: we call it the {\it singular} companion of $L$ and denote it ${\rm sing}(L)$.
The following proposition holds by Proposition 9.4 and Remark 9.5 {\it loc. cit.} 

\begin{prop} 
\label{prop:singcomp} 
Assume $n \equiv 4 \bmod 8$ and denote by $\mathcal{A}$ the 
set of isometry classes of exceptional $L$ in $\mathcal{L}_n$ with ${\rm r}_1(L)=0$, 
and by $\mathcal{B}$ the set of isometry classes of non exceptional $L$ in $\mathcal{L}_n$ with ${\rm r}_1(L)\neq 0$.
Then we have a natural bijection $\mathcal{A} \rightarrow \mathcal{B}, [L] \mapsto [{\rm sing}(L)]$. 
Moreover, for $L \in \mathcal{A}$ and $L':={\rm sing}(L)$ we have 
$|{\rm Exc}(L)|={\rm r}_1(L')$, ${\rm R}_2(L)={\rm R}_2(L')$ and $|{\rm O}(L')|\,=\,2 |{\rm O}(L)|$.
\end{prop}

In particular, we always have $|{\rm Exc}(L)|=2m$ with $0 \leq m\leq n$ for $L \in \mathcal{L}_n$.

\section{The invariant ${\rm BV}$}
\label{sect:inv}

Let $G$ be an arbitrary finite graph with set of vertices $V$. Let $A$ be the adjacency matrix of $G$, 
a $V \times V$ matrix.\footnote{In the application below, $G$ will be undirected
and for all $v,v' \in V$ we will have at most one edge between $v$ and $v'$.} Consider the square 
$S:=A^2$ of $A$, say $S=(s_{u,v})_{(u,v) \in V \times V}$. For $v \in V$, we define $${\rm C}(G;v)=\{\{ \,s_{u,v} \, |\, u \in V\}\}$$ as the multiset of entries of the column $v$ of the matrix $S$. In other words, ${\rm C}(G;v)$ is the multiset of {\it numbers of length $2$ paths in $G$ starting at $v$ and ending at another given vertex}. 

\begin{definition} 
\label{def:BVG}
Let $G$ be a graph with finite set of vertices $V$.
For $v \in V$, we define ${\rm BV}(G)$ as the multiset $\{\{ {\rm C}(G;v)\, |\, v \in V\}\}$. 
This is a multiset of multisets of integers, and it is an invariant of the isomorphism class of $G$.
\end{definition}

Assume now $L$ is an Euclidean integral lattice. 
For $i\geq 0$, we view the finite set ${\rm R}_{\leq i}(L)=\{ v \in L\, \, |\, \, v.v \leq i\}$ 
as a metric set (or better, as an {\it Euclidean set} as in \cite[\S 4.1]{cheuni}). 
We denote by\footnote{This is a variant of the graph also denoted by ${\rm G}(L)$
in \cite{cheuni} \S 4.} ${\rm G}(L)$ the undirected graph whose vertices are the nonzero 
pairs $\{\pm v\}$ with $v \in {\rm R}_{\leq 3}(L)$, and with 
an arrow between $\{ \pm v\}$ and $\{\pm w\}$ if, and only if, $v.w \equiv 1 \bmod 2$.

\begin{definition} 
\label{def:BVL}
If $L$ is an Euclidean integral lattice. 
We define ${\rm BV}(L)$ as the multiset of multisets of integers ${\rm BV}({\rm G}(L))$.
\end{definition}

By construction, ${\rm BV}(L)$ is an invariant of the isometry class of ${\rm R}_{\leq 3}(L)$, hence of $L$.
We choose the notations ${\rm BV}$ for {\it Bacher-Venkov}, as this definition is inspired from an invariant defined in \cite{bachervenkov}.

\begin{remark} {\rm (The Bacher-Venkov polynomials)} Assume we are in the special case ${\rm R}_{\leq 2}(L)=\emptyset$. If we choose distinct $\pm u$ and $\pm v$ in ${\rm R}_3(L)$, 
we have either $u.v=0$ or $u.v = \pm 1$, as $u.v=2$ implies $u-v \in {\rm R}_2(L)$. 
For $v\in {\rm R}_3(L)$, Bacher and Venkov define in {\rm \cite[p.15]{bachervenkov}} the polynomial 
${\rm m}_v(x)$ as the sum, over all $w \in {\rm R}_3(L)$ with $w.v=1$, of $x^{{\rm n}(v,w)}$, with ${\rm n}(v,w)=| \{ \xi \in {\rm R}_3(L)\, |\, v.\xi=w.\xi=-1\}|$. If, in this definition of ${\rm m}_v(x)$, we rather sum over all $w \in {\rm R}_3(L)$, the obtained polynomial contains the same information as the multiset ${\rm C}({\rm G}(L); v)$.
\end{remark}

We now give a few information about the graph ${\rm G}(L)$ for $L \in \mathcal{L}_{n}$
with ${\rm r}_1(L)=0$ and $n=28,29$. Its number of vertices is $\frac{1}{2}({\rm r}_2(L)+{\rm r}_3(L))$ and we know from \cite{king} the possibilities for ${\rm R}_2(L)$. \ps\ps

(a) ({\rm \bf The case $n=28$)} A theta series computation shows ${\rm r}_3(L)\,=\,2240\,+\, 8 \,{\rm r}_2(L) -\,256\,|{\rm Exc}(L)|$ (see \cite[\S 4]{bachervenkov} and \cite[Prop. 4.8]{cheuni}). 
Here, $|{\rm Exc}(L)|$ is the number of characteristic vectors of norm $4$ of $L$; it satisfies $0 \leq |{\rm Exc}(L)| \leq 56$.
We refer to Proposition~\ref{prop:singcomp} for this inequality, and to Sect. \ref{sect:notation} for the notion of singular companion used in the following example.

\begin{example}
\label{ex:r3zero}
{\rm (The case ${\rm r}_3=0$)} 
{\it There are exactly $15$ classes of lattices $L \in \mathcal{L}_{28}$ with ${\rm r}_1(L)={\rm r}_3(L)=0$.
Indeed, the formula above for ${\rm r}_3(L)$ shows $|{\rm Exc}(L)|\geq 10$ {\rm (}in particular, $L$ is exceptional{\rm )}. 
Let $L'$ denote the singular companion of $L$, say $L' \simeq {\rm I}_m \perp U$ with $U$ unimodular of rank $28-m$ with ${\rm r}_1(U)=0$ and $2m=|{\rm Exc}(L)|$. We have ${\rm r}_3(L)=0$ if and only if ${\rm r}_2(U)\,=2\,m\,(33-m)-280$. 
We easily conclude using the classification of unimodular lattices of rank $\leq 23$. The extreme cases are $L'={\rm I}_{28}$ {\rm ($m=28$)} and $L' = {\rm I}_5 \perp U$ with $U$ the short Leech lattice {\rm ($m=5$)}. These $15$ lattices have different root systems.}
\end{example}

As an indication, it follows from our final computations that ${\rm G}(L)$ has between $20$ and $3388$ vertices, and $1318$ in average, for $L \in \mathcal{L}_{28}$ with ${\rm r}_1(L)=0$. Better, for $97$\% of these lattices this number lies in $[1000, 1600[$:

\tabcolsep=1.8pt
\begin{table}[H]
{\scriptsize 
\renewcommand{\arraystretch}{1.8} \medskip
\begin{center}
\begin{tabular}{c|c|c|c|c|c|c|c|c|c|c|c|c|c|c|c|c|c}
$N$ & $0$ & $200$ & $400$ & $600$ & $800$ & $1000$ & $1200$ & $1400$ & $1600$ & $1800$ & $2000$ & $2200$ & $2400$ & $2600$ & $2800$ & $3000$ & $3200$ \\
$\sharp$ & $12$ & $24$ & $66$ & $273$ & $1801$ & $22889$ & $269512$ & $53898$ & $6805$ & $1269$ & $297$ & $84$ & $51$ & $12$ & $5$ & $1$ & $4$
\end{tabular} 
\end{center}
\caption{The number $\sharp$ of $L \in {\rm X}_{28}$ with ${\rm r}_1(L)=0$ such that the number of vertices of ${\rm G}(L)$ lies in the interval $[N, N+200[$.}
\label{tab:sizeR3}
}
\end{table}	

The average number of edges of those graphs is $\simeq 470\,000$. 
Another interesting property,\footnote{It may be possible to explain part of these observations by using harmonic theta series arguments as in \cite{venkov}, but we shall not pursue this (unnecessary) direction here.} 
is that ${\rm R}_{\leq 3}(L)$ does generate $L$ over $\Z$ for most lattices, 
as indicated by Table~\ref{tab:indexR3}. For example, the only lattice such that ${\rm R}_{\leq 3}(L)$ does not generate a finite index subgroup of $L$ is ``the'' exceptional lattice with companion ${\rm I}_5 \perp U$, with $U \in {\rm X}_{23}^\emptyset$ the Odd Leech lattice. 

\tabcolsep=4pt
\begin{table}[H]
{\scriptsize 
\renewcommand{\arraystretch}{1.8} \medskip
\begin{center}
\begin{tabular}{c|c|c|c|c|c|c|c|c|c|c|c|c|c|c|c|c}
 $d$  & $1$ & $2$ & $3$ &  $4$ & $5$ & $6$ & $7$ & $8$ & $9$ & $10$ & $12$ & $13$ & $16$ & $20$ & $24$  & $25$  \\
 $\sharp $  & $356462$ & $364$ & $16$ & $58$ & $5$ & $8$ & $3$ & $18$ & $7$ & $1$ & $9$  & $2$ & $9$ & $5$ & $1$  & $2$ \\
 \hline
 $d$  & $27$ & $32$ & $36$ & $49$ & $64$ & $72$ & $125$ & $128$ & $243$ & $256$ & $729$ & $2048$ & $4096$ & $8192$ & $\infty$ & other  \\
 $\sharp $  & $2$ & $8$ & $1$ &  $2$ &  $3$ & $4$ & $2$ & $3$ & $1$ & $2$ & $1$ & $1$ & $1$ & $1$ & $1$ & $0$
\end{tabular} 
\end{center}
\caption{The number $\sharp$ of $L \in {\rm X}_{28}$ with ${\rm r}_1(L)=0$ such that ${\rm R}_{\leq 3}(L)$ generates a sublattice of index $d$ in $L$.}
\label{tab:indexR3}
}
\end{table}

(b) ({\rm \bf The case $n=29$)} For $L$ in $\mathcal{L}_{29}$ with ${\rm r}_1(L)={\rm r}_2(L)=0$, 
a theta series computation shows ${\rm r}_3(L)\,=\,1856\,-\,128\, |{\rm Exc}(L)|$, where ${\rm Exc}(L)$ is the number of characteristic vectors of norm $5$ in $L$. It is a fact, that we shall not explain here, that we always have $|{\rm Exc}(L)| \leq 2$. So we have either $\frac{1}{2}{\rm r}_3(L)=928$ in the non exceptional case, and $\frac{1}{2}{\rm r}_3(L)=800$ otherwise. The situation here is thus much more uniform than in the case $n=28$. Also, the graph ${\rm G}(L)$ 
turns out to have $259\,840$ edges in the non exceptional case, and $198\,400$ otherwise, and in all cases ${\rm R}_3(L)$ does generate $L$ over $\Z$.\ps

\ps	
{\bf Some computational aspects of ${\rm BV}$.} 
It is straightforward in principle to compute ${\rm BV}(L)$ from a given Gram matrix $M$ of the integral lattice $L$.
Indeed, we can use the Fincke-Pohst algorithm to find a column matrix $V$ 
whose raws are the $\pm v$ with $v \in L$ such that $0<v.v \leq 3$. 
The adjacency matrix $A$ of ${\rm G}(L)$, which is symmetric and with coefficients $0$ or $1$,
is determined by $A \equiv {}^{\rm t}V \,M \,V \bmod 2$ (a fast computation). 
More lengthy is then the computation of $S=A^2$, since the size of $A$ is typically a thousand or more. 
To save time, this squaring is computed only modulo some large enough prime (we used $1009$), 
and is performed using only single-word arithmetic.\footnote{
While GP normally use multiprecision integer arithmetic, 
we used single word arithmetic (using \texttt{t\textunderscore VECSMALL})
whenever appropriate in all of our programs.} 
In practice, the slightly weaker invariant than ${\rm BV}(L)$ consisting of the {\it set} (rather than multiset) of the multisets  
${\rm C}({\rm G}(L); v)$ (with $v \in V$) proved equally strong, and this what we implemented in our computations.
Finally, we apply a hash function to the resulting invariant to get a $64$bit identifier. This allows
not only to quickly discard lattices with already known identifiers, but also to store the list of known identifiers in a compact
way. The choice of $1009$ and of the hash function here is arbitrary, all that matter is
that the resulting invariant is fine enough to distinguish all the lattices we consider, which can be only determined a posteriori.

\ps\ps
Here are a few CPU time information in our range:\ps

-- For each of the $346\,299$ lattices $L \in {\rm X}_{28}$ with ${\rm r}_1(L)=0$ and such that the number of vertices of ${\rm G}(L)$ is in $[1000, 1600[$, the average CPU time to compute ${\rm BV}(L)$ is about \texttt{1.5 s}. For the worst (and actually irrelevant!) case with $3388$ vertices, the CPU time is approximately $\texttt{35 s}$. \ps

-- For each of the $10092$ lattices $L \in {\rm X}_{29}^\emptyset$, the average CPU time to compute ${\rm BV}(L)$ is about \texttt{1.2 s}. \ps


\section{Finding $\Z$-basis made of short vectors}
\label{sect:Zbasisshort}

An important ingredient in our computation is the Plesken-Souvignier algorithm \cite{pleskensouvignier}, 
which computes $|{\rm O}(L)|$ from the Gram matrix of a given $\Z$-basis $e=(e_1,\dots,e_n)$ of $L$. For this algorithm to be efficient (and not too memory consuming), it is crucial to have ${\rm m}(e):={\rm Max}\, \{ e_i.e_i\, \, |\, \, i=1,\dots,n\}$ as small as possible, and also highly desirable to have $|\{ 1 \leq i \leq n\,\, |\, e_i.e_i= {\rm m}(e)\}|$ small. The {\rm LLL} algorithm, although very fast and useful, does not provide in general a $\Z$-basis of $L$ which is good enough in these respects.

\begin{example} 
\label{ex:LLLU29} For our $10092$ elements in ${\rm X}_{29}^\emptyset$,
the ${\rm LLL}$ algorithm produces for $6570$ of them a basis $e$ with ${\rm m}(e)=4$, for $3512$ a basis with ${\rm m}(e)=5$, and in the $10$ remaining cases a basis $e$ with ${\rm m}(e)=6$. As already said, those lattices are actually generated over $\Z$ by their ${\rm R}_{\leq 3}$, hence may (and actually do) have a $\Z$-basis $e$ with $m(e)=3$. 
\end{example}

In order to search for better lattice bases, we used the following simple probabilistic algorithm. 
Its main function $\texttt{reduce}$ takes as an input an integral Euclidean lattice $L$ of rank $d$, 
an integer $b$, and another integer $t$ ("number of tries"). 
It returns either $0$ (failure) or a basis $e$ of $L$ with ${\rm m}(e)=b$. 
\ps

1. Compute the set $S$ of vectors $\pm v$ in $L$ with $0<v.v\leq b$. \par
2. If $S$ does not generate $L$ over $\Z$, return $0$.\par
3. Compute the set $R$ of vectors $\pm v$ in $L$ with $0<v.v \leq b-1$.\par
4.  Compute the rank\footnote{Thus step is especially useful in the case $b=3$, but may be ignored in situations where we know that $R$ is big (in which case we set $k_0=1$).} $r$ of $R$ and set $k_0={\rm Max}(1,d-r)$. \par
5.  For $k$ from $k_0$ to $d$, for $i$ from $1$ to $t$, do\par
\indent \hspace{1cm} 5a. Choose $k$ vectors $e_1,\dots,e_k$ randomly in $S$,\par
\indent \hspace{1cm} 5b. Choose $d-k$ vectors $e_{k+1},\dots,e_d$ randomly in $R$,\par
\indent \hspace{1cm} 5c. If $e_1,\dots,e_d$ generates $L$ over $\Z$, return $e_1,\dots,e_d$.\par
6. Return $0$.\ps

For a given integral lattice $L$, we start with a $\Z$-basis $e$ given by the {\rm LLL} algorithm, 
choose some $t$, and then we apply $\texttt{reduce}(L,b,t)$ successively to $b=1,\dots,{\rm m}(e)-1$ until it returns some $\Z$-basis of $L$. If it fails, it means it did not beat the initial basis $e$ given by {\rm LLL}. 


\begin{remark}
\label{rem:detreduce}
{\rm Step 2 and 5c of the function \texttt{reduce} amounts to checking that some determinant is $\pm 1$.
To save time we may first check that this holds modulo $2$ or other primes.}
\end{remark}

\begin{example} 
\label{ex:reduceU29}
{\rm For each lattice $L$ in ${\rm X}_{29}^\emptyset$, the $\texttt{reduce}$ algorithm does find a $\Z$-basis consisting of norm $3$ vectors of $L$ in about \texttt{93 ms} (using $t=1000$ is usually enough, and computing the set of norm $3$ vectors already takes about \texttt{30 ms}). Using those bases, the average time to determine $|{\rm O}(L)|$ using the Plesken-Souvignier algorithm\footnote{We used the $\texttt{PARI/GP}$ implementation $\texttt{qfauto}$ of Souvignier's code, with flag $[0,2]$.}
 is about \texttt{1.24 s}.}
\end{example}

	We now discuss the $357003$ rank $28$ unimodular lattices with no norm $1$ vectors.
As already explained in Table~\ref{tab:indexR3}, $356462$ of them are generated over $\Z$ 
by their ${\rm R}_{\leq 3}$, and we did find a $\Z$-basis of norm $\leq 3$ vectors in all cases, in about $
\texttt{158 ms}$. For the $541$ remaining lattices, $519$ are generated over $\Z$ by their ${\rm R}_{\leq 4}$ and we also found in all cases a $\Z$-basis consisting of norm $\leq 4$ vectors using $\texttt{reduce}(L,4,10000)$.\ps
	The remaining $22$ lattices are atypical, but not especially mysterious. They all have a root system of rank $28$, except for two of them, whose root systems are respectively  $2{\bf A}_7\,{\bf D}_{13}$ (rank $27$) and ${\bf D}_5$ (which has only rank $5$, but this lattice is the exceptional lattice with companion ${\rm I}_{5} \perp {\rm Odd Leech}$ discussed in Example~\ref{ex:r3zero}). For instance, for many of these lattices, the sublattice generated by the root system ${\rm R}_2(L)$ is isometric to ${\rm D}_{m_1} \perp  {\rm D}_{m_2} \perp \cdots \perp {\rm D}_{m_s}$ with $m_1+\dots+m_s=28$, where ${\rm D}_m \subset {\rm I}_m$ is the standard root lattice of type ${\bf D}_m$, namely  ${\rm D}_m={\rm M}_2(1^m)$. But for $m\geq 1$, the dual lattice ${\rm D}_m^\sharp$ writes
$${\rm D}_m^\sharp = {\rm D}_m \coprod (\varepsilon_m+{\rm D}_m) \coprod (\eta^+_m+{\rm D}_m) \coprod (\eta^-_m+{\rm D}_m) $$
with $\varepsilon_m=(0,\dots,0,1)$ and $\eta^{\pm}_m=\frac{1}{2}(1,\dots,1,\pm 1)$. As is well-known, the minimum of $x \mapsto x.x$ on $\eta^{\pm}_m+{\rm D}_m$ is $m/4$ (resp. $1$ on $\varepsilon_m+{\rm D}_m$). Assuming ${\rm r}_1(L)=0$, it follows that the small norm vectors of those $L$ tend to generate proper sublattices of $L$. 
The two most striking cases are the following:
	
\begin{example} 
\label{ex:extremeD28}
\begin{itemize}
\item[(a)] There is a unique $L \in {\rm X}_{28}$ with ${\rm r}_1(L)=0$ and root system ${\bf D}_{28}$. We may define it as
$L = {\rm D}_{28} \coprod C$ with $C=\eta^+_{28} + {\rm D}_{28}$. 
As ${\rm min}_{x \in C}\,  x.x\,=\,28/4\,=7$, it follows that $L$ is generated over $\Z$ by the set ${\rm R}_{\leq 7}(L)$, whereas ${\rm R}_{\leq 6}(L)$ generates the index $2$ subgroup ${\rm D}_{28}$. Note that ${\rm R}_{\leq 7}(L)$ is huge: it has $158\,736\,881$ vectors! \ps
\item[(b)] There is a unique $L \in {\rm X}_{28}$ with ${\rm r}_1(L)=0$ and root system ${\bf D}_{8}\,{\bf D}_{20}$. 
Set $Q={\rm D}_{8} \perp {\rm D}_{20}$. We may define $L \subset Q^\sharp$ as the inverse image, under the natural map $Q^\sharp \rightarrow Q^{\sharp}/Q = {\rm D}_8^\sharp/{\rm D}_8 \oplus {\rm D}_{20}^\sharp/{\rm D}_{20}$, of the bilinear Lagrangian 
$\{ 0,\, \varepsilon_{8}\,+\,\eta^{+}_{20}, \,\eta^+_8\,+\,\varepsilon_{20},\, \eta^{-}_8\,+\,\eta^-_{20}\}$ {\rm (}see~{\rm \cite[\S 2]{cheuni}}{\rm )}. It follows that $L$ is a union of $4$ cosets of $Q$ with respective minimum norm 
$0$, $1+20/4=6$, $8/4+1=3$ and $8/4+20/4=7$. This shows that $L$ is generated by ${\rm R}_{\leq 6}$ (with $21\,827\,953$ vectors!) but not by ${\rm R}_{\leq 5}$.
\end{itemize}
\end{example}
	
These are the worst cases, since the $20$ other lattices are generated over $\Z$ by their ${\rm R}_{\leq 5}$ and we did find for all of them a $\Z$-basis consisting of norm $\leq 5$ vectors using $\texttt{reduce}(L,5,10000)$. 
We refer to~\cite{allcheweb28} for a list of Gram matrices that we found using these methods. \ps

We finally discuss the computation of the order of the isometry group ${\rm O}(L)$ for $L$ in ${\rm X}_{28}$.
It is enough to compute the order of the {\it reduced isometry group} ${\rm O}(L)/{\rm W}(L)$ of $L$ (see Sect.~\ref{sect:notation}).
As already observed in~\cite{chniemeier}, this is much more efficient, and can easily be done using features of the Plesken-Souvignier algorithm: see~\cite[\S 4.2]{cheuni}. Using the Gram matrices above for the lattices in ${\rm X}_{28}$, this is quite fast.   
For instance, for $99.8 \%$ of the $356462$ lattices generated over $\Z$ by their ${\rm R}_{\leq 3}$, the average time to determine their reduced isometry group is about \texttt{0.3 s}. Note that in the two extreme cases discussed in Example~\ref{ex:extremeD28}, the reduced isometry group is clearly trivial, so there is nothing to compute.


\section{The Biased Neighbor Enumeration algorithm}
\label{sect:thebiasedalgo}
 In this section, we give in an explicit form the general algorithm informally described in \cite[\S 1.10]{cheuni}. 
 Its goal is to produce large quantities of cyclic $d$-neighbors of ${\rm I}_n$ having a given root system, by imposing in the enumeration of $d$-isotropic lines a suitably chosen {\it visible root system} in the sense of {\it loc. cit.} \S 5. 
 We start with a few notation and terminology. \ps\ps

For $x \in \R^n$ and $j \in \R$, we denote by ${\rm m}_j(x) \in \{0,1,\dots,n\}$ the number of coordinates of $x$ which are equal to $j$. A pair $(d,x)$, with $d\geq 1$ an integer and $x \in \Z^n$, will be said {\it normalized} if we have 
$1=x_1 \leq x_2 \leq \cdots \leq x_n \leq d/2$, as well as ${\rm m}_1(x) \geq {\rm m}_j(x)$ for all $1 \leq j <d/2$.
We call ${\rm m}_1(x)$ the {\it index} of $(d,x)$, and ${\rm m}_{d/2}(x)$ the {\it end}\footnote{Obviously,  the end of $(d,x)$ is nonzero only for $d$ even and $x_n=d/2$.} of $(d,x)$. Finally, the {\it type} of $x$, or of $(d,x)$, is the partition of the integer $n$ defined by the nonzero integers ${\rm m}_j(x)$, $j \in \Z$.  

\begin{example} 
\label{ex:exemplexd}
Let $x=(\scalebox{.7}{{\rm 1, 1, 1, 2, 3, 4, 4, 5, 6, 7, 8, 8, 9, 9, 10, 10, 10, 11, 11}}) \in \Z^{19}$ and $d=22$. 
Then $(d,x)$ is normalized, of type {\small $3+3+2+2+2+2+1+1+1+1+1$}, index $3$ and end $2$. 
We also use the notation $3^2\,2^4\,1^5$ for such a partition.
\end{example}

The basis reason for this terminology the presence of a large isometry group of the lattice ${\rm I}_n=\Z^n$, namely the $2^n n!$ permutations and sign changes of coordinates, which has the following immediate consequence:

\begin{fact}
\label{fact:orbits} Let $d\geq1$ and $x \in \Z^n$. Assume $x_i \not\equiv 0 \bmod d$ for all $1 \leq i \leq n$, and that there exists\footnote{Note that this condition is automatically satisfied if $d$ is prime.} $1\leq j \leq n$ such that $a:=x_j$ is prime to $d$ and satisfies ${\rm m}_a(x)\geq {\rm m}_{i}(x)$ for all $i \in \Z$. Then there is $\sigma \in {\rm O}({\rm I}_n)$ and $y \in \Z^n$ such that $(d,y)$ is normalized and satisfies $\sigma(x) \equiv a y \bmod d \Z^n$.
\end{fact}


We now describe the Biased Neighbor Enumeration algorithm, later refer to as \texttt{BNE}.
It takes as inputs: an integer $n \geq 1$, an isomorphism class of root system ${\rm R}$, 
the reduced mass $\rho$ of ${\rm X}_n^R$, an integer $0\leq e<n$ and 
integer partition $n-e=\sum_{i=1}^s n_i$, with $n_1 \geq n_2 \geq \cdots \geq n_s\geq 1$. 
If \texttt{BNE} terminates, it returns a list of $(d,x,\epsilon,\mu,\beta)$ where 
the lattices ${\rm N}_d(x; \epsilon)$ are representatives of ${\rm X}_n^R$, and where ${\rm N}_d(x;\epsilon)$ has reduced mass $\mu$ and {\rm BV} invariant $\beta$; it also proves that all the elements in ${\rm X}_n^R$ have distinct ${\rm BV}$ invariants. \ps \ps
{\small 
1. Set $d=2$ and define empty lists $inv$ and $lat$.\par
2. Make the list $L$ of all $x \in \Z^n$ with $(d,x)$ normalized of type $n_1+\dots+n_s+e$, index $n_1$ and end $e$.\par
3. Only keep in $L$ those $x$ such that $x$ is $d$-isotropic. \par
4. For each $x \in L$ and\footnote{We restrict to $\epsilon=0$ for $d$ odd or $e>0$.} $\epsilon \in \{0,1\}$, compute the cyclic $d$-neighbors ${\rm N}_d(x; \epsilon)$ of ${\rm I}_{n}$.\par
5. Replace $L$ with the list of $(x,\epsilon) \in L \times \{0,1\}$, satisfying both ${\rm r}_1({\rm N}_d(x;\epsilon))=0$ and ${\rm R}_2({\rm N}_d(x;\epsilon)) \simeq R$.\par
6.  Compute the (hashed) ${\rm BV}$ invariants $\beta$ of all ${\rm N}_d(x;\epsilon)$ with $(x,\epsilon) \in L$.\par
7. For each $(x,\epsilon) \in L$, if the BV invariant $\beta$ of ${\rm N}_d(x;\epsilon)$ is not in $inv$ do:\par
${}^{}$ \indent  7a. Add $\beta$ to the list $inv$,\par
${}^{}$ \indent  7b. Compute the reduced mass $\mu$ of ${\rm N}_d(x;\epsilon)$,\par
${}^{}$ \indent  7c. Add $(d,x,\epsilon,\beta,\mu)$ to the list $lat$,\par
${}^{}$ \indent  7d. Set $\rho \leftarrow \rho-\mu$. If $\rho=0$, return $lat$.\par
8. $d \leftarrow d+1$ and go back to Step 2.\ps
}
\ps\ps

The main idea of this ``Coupon Collector'' algorithm is explained in details in~\cite{cheuni} \S 1.10 and Sect. 5.
We will briefly review it below, including the (key) role of the $n_i$ and $e$, but we first state an important criterion for {\rm \texttt{BNE}} to terminate.
For an ADE root system $S$, and an integral lattice $L$, we denote by ${\rm emb}(S,L)$ 
the number of isometric embeddings\footnote{See the general notation at the end of Sect. \S~\ref{sect:intro}.} 
${\rm Q}(S) \rightarrow L$ with saturated image. 
Following \S 5.9 {\it loc. cit.}\,, we say that a pair $(R,S)$ of ADE root systems is {\it safe} 
if for any integral lattices $L$ with ${\rm r}_1(L)=0$ and ${\rm R}_2(L) \simeq R$ we have ${\rm emb}(S,L) \neq 0$.
We attach to $e$ and the integer partition $n=n_1+\cdots+n_s+e$ the 
root system $V:= {\bf A}_{n_1-1}{\bf A}_{n_2-1} \cdots {\bf A}_{n_s-1}{\bf D}_e$.
By Theorem E of \cite{cheuni}, a special case of the results of \cite{chestat}, we have:

\begin{thm} 
\label{thm:algoterm}
Assume $(R,V)$ is safe and that different classes in ${\rm X}_n^R$ have different ${\rm BV}$ invariants.
Then the algorithm {\rm \texttt{BNE}} terminates.
\end{thm}

Let us explain this theorem. 
Recall that the {\it visible root system} of a $d$-neighbor $N$ of ${\rm I}_n$ is defined as 
$R^{\rm v}:={\rm R}_2(N) \cap {\rm I}_n$. We refer to {\it loc. cit.} Sect. 5 for a study of its specific properties.
Let us simply say here that $R^{\rm v}$ tends to generate a saturated sub lattice in $N$ (this is holds at least when $d$ is a large prime). By Fact~\ref{fact:orbits}, Steps 2 and 3 enumerate 
representatives $(d,x,\epsilon)$ for almost all ${\rm O}({\rm I}_n)$-orbits of 
the cyclic $d$-neighbors ${\rm  N}_d(x;\epsilon)$ of ${\rm I}_n$ 
whose {\it visible root system} is isomorphic to $V$.
In the case $d$ is prime, all orbits are considered in these two steps.
By the aforementioned Theorem E {\it loc. cit.}, for any given $L \in {\rm X}_n$, 
the proportion of $d$-neighbors of ${\rm I}_n$ having a given visible root system $V$, and which are isomorphic to $L$, 
tends to (an absolute constant times) ${\rm emb}(V,L)/|{\rm O}(L)|$ when the prime $d$ goes to $\infty$.
This number is non zero when $(R,V)$ is safe, which shows that {\rm \texttt{BNE}} terminates, hence Theorem~\ref{thm:algoterm}. We stress however that {\it there is no known upper bound on the integer $d$ such that 
each $L \in \mathcal{L}_n$ is isometric to a $d$-neighbor of ${\rm I}_n$}. 
This method for searching for unimodular lattices is probabilistic and a form of (non-uniform) {\it coupon collector problem}: see \cite[\S 1.12]{cheuni} for a discussion along these lines.

\begin{remark} {\rm (Choice of $V$)} 
\label{rem:choiceV}The root system $R$ being given, it is always possible to choose $V$ such that $(R,V)$ is safe.
For instance, choosing $V=\emptyset$, {\it i.e.} $e=0$ and all $n_i$ equal to $1$, trivially works, although it is usually inefficient if $R$ is large. Indeed, the algorithm will most likely find lattices with smaller root systems than $R$, as we have ${\rm emb}(V,L)=1$ for all $L$. In practice, the game is rather to choose for $R^{\rm v}$ a sub-root system of $R$ which is as large as possible, given the general constraints for visible root systems: see Sect. 5 {\it loc. cit.} for a discussion of good and possible choices. See also the next section for examples.
\end{remark} \ps 

	We finally discuss certain steps or features of \texttt{BNE}.

\begin{remark} 
\label{rem:linesversusvectors}
{\rm (Step $2'$:  lines versus vectors, and redundancy)} 
{\rm
Let $x,x' \in \Z^n$ be two $d$-isotropic vectors, say with $(d,x)$ and $(d,x')$ normalized 
of same type $n_1+\dots+n_s$. 
Denote by $\ell$ and $\ell' \subset {\rm I}_n \otimes \Z/d$ the $\Z/d$-line they generate.
If $\ell$ and $\ell'$ are in the same ${\rm O}({\rm I}_n)$-orbit, then the $d$-neighbors defined by 
$(d,x)$ and $(d,x')$ are naturally isomorphic, creating unwanted redundancy in our algorithm. 
For a given $x$, then there are exactly $f$ vectors $x'$ generating an equivalent line in this sense,
with $f:=|\,\{\,1\leq i \leq d/2\,\,|\,\,{\rm m}_i(x)={\rm m}_1(x)\,\,\& \,\,{\rm gcd}(i,d)=1\}\,|$. 
This is especially large in the case $n_i=1$ for all $i$, and $d$ is prime, 
for which we have $f=n$. We did not find any clever way to directly select a single element among those $f$ ones, but did instead the following right after Step 2 in the algorithm: fix some total ordering $\prec$ on $\Z^n$, and for each $x \in L$, 
compute all $x'$ equivalent to $x$ in the sense above, and only keep $x$ in $L$ if it is the biggest of them for $\prec$. }
\end{remark}

\begin{remark} {\rm (Step 5) See {\it e.g.} Remark 4.4. in \cite{cheuni} for an algorithm to compute root systems.}
\end{remark}

\begin{remark} {\rm (Step 7b)
\label{rem:step7b}
For this step, of course, we apply the ideas exposed in Sect.~\ref{sect:Zbasisshort}: 
for each $L:={\rm N}_d(x,\epsilon)$, we first search for a $\Z$-basis of $L$ made of small norm vectors with \texttt{reduce} beating \texttt{LLL}, and then we compute the order of the reduced isometry group ${\rm O}(L)/{\rm W}(L)$ by giving this basis to the Plesken-Souvignier algorithm as explained in \cite[\S 4.2]{cheuni}. 
}
\end{remark}

\begin{remark} {\rm (Parallelization)
\label{rem:parallel}
Each of Steps $3$ to $6$ (as well as Step $7b$) is straightforward to parallelize in practice. 
For memory reasons, we usually modify Step $2$ by limiting the size of the list $L$ 
and go to Step 3 when this limit is reached. If so, at the end of Step 7, 
we go back to Step 2 and pursue the enumeration of the remaining $x$ 
until all of them have been considered, before going to Step 8. This way to present the algorithm 
is especially suited to find the large (and most difficult) ${\rm X}_n^R$. 
In those cases, most of the CPU time of \texttt{BNE} is spent on Step 6.
This way of limiting the size of $L$ also allows of course to use 'early abort' strategy and cut search time 
when a single lattice is missing (which eventually always happens!). 
}
\end{remark}

	Of course, the algorithm is interesting even if it does not terminate, since it usually finds many lattices if not all.
Also, in most cases we also do not really start at $d=2$ but at $d=2s+1$, otherwise there is clearly no normalized pair $(d,x)$ of type $n_1+\dots+n_s+e$. For some purposes, we may also restrict the enumeration of $(d,x)$ to $d$ in certain congruence classes modulo some integers.
For instance in some hard cases, we sometimes restrict to $d$ odd to optimize further $2$-neighbor computations. 
We must also have $d$ even in the case $e\neq 0$. 
In principle, we could take $d$ large from the beginning, but this goes against our (guilty) wish to find neighbor forms of smallest possible ``farness'' for our lattices (that is of the form ${\rm N}_d(x;\epsilon)$ with $d$ as small as possible).\ps

\begin{remark} {\rm (The non-biased $\texttt{NE}$ algorithm) The simpler, non biased, variant $\texttt{NE}$ of $\texttt{BNE}$, 
is the case where we enumerate in Step 2 the pairs $(d,x)$ of all possible types, and seek for all root systems at the same time. The only input of $\texttt{NE}$ is the integer $n$ and the dictionnary $M : R \mapsto $ reduced mass of ${\rm X}_n(R)$, and it returns ${\rm X}_n'$ (and theoretically terminates). The only difference is that in Step 5 we only keep $(x,\epsilon)$ if the root system of ${\rm N}_d(x;\epsilon)$ still belongs to $M$, and in Step 7d we delete the root system $R$ from $M$ if the remaining reduced mass of ${\rm X}_n(R)$ is zero after this step. 
This is how we actually started our search of ${\rm X}_{28}$. Although $\texttt{NE}$ allows to fills quite many ${\rm X}_{28}^R$ having a small quantities of lattices, it becomes very lengthy and inefficient when $d$ grows as explained in \cite{cheuni} \S 1, and it is absolutely curcial to use $\texttt{BNE}$ instead.}
\end{remark}

\section{Some examples}
\label{sect:examples}
${}^{}$ \par

We start with a few simple examples for which the $\texttt{BNE}$ algorithm directly works,
and then provide a few more complicated ones. 
Note that there are usually many ways to find a given lattice, so some of the lattices described below may appear in a different neighbor form in our lists \cite{allcheweb28} and \cite{allcheweb29}.  \ps

\subsection{The root system $7{\bf A}_1\,3{\bf A}_2\,{\bf A}_7$ in dimension $28$}

${}^{}$ \indent For this root system $R$, we know from King that the reduced mass of ${\rm X}_{28}^R$ is $5/24$ (whereas its mass is 
$|{\rm W}(R)|=1\,114\,767\,360$ times smaller), so we expect ${\rm X}_{28}^R$ to contain only very few lattices.
We choose the visible root system $V:=6{\bf A}_1\,2{\bf A}_2\,{\bf A}_7$, that is $e=0$ and the integer partition 
$8\,3^2\,2^6\,1^2$ of $28$. We have $s=11$ so we may start at $d=23$. 
The $\texttt{BNE}$ algorithm terminates at $d=27$ and returns the $3$ lattices in ${\rm X}_{28}^R$, with reduced masses {\small $1/12$}, {\small $1/12$} and {\small $1/24$}, after about $7$ minutes of CPU time. \ps

It will be convenient to use the following kind of tables to give some details about the intermediate steps of such calculations. The second column gives, for the $d$ in the first column, the number $\sharp \texttt{iso}$ of $d$-isotropic lines found after Steps $2$ and $3$, incorporating Step $2'$ explained in Remark~\ref{rem:linesversusvectors}. The third column gives the size $\sharp\texttt{found}$ of the list $L$ after Step $5$. The better we have chosen $V$, the higher the ratio $\sharp \texttt{iso}/\sharp\texttt{found}$ is. The fourth column gives the 
number of new lattices found in $L$, and the last column, the remaining reduced mass $\rho$ of ${\rm X}_{n}^R$ after Step $7$ (when $0$, the algorithm terminates).

\tabcolsep=7pt
\begin{table}[H]
{\scriptsize 
\renewcommand{\arraystretch}{1.2} \medskip
\begin{center}
\begin{tabular}{c|c|c|c|c}
$d$ & $\sharp \texttt{iso}$ & $\sharp\texttt{found}$ & $\sharp \texttt{new\textunderscore lat}$ & $\texttt{rem\textunderscore red\textunderscore mass}$  \\
\hline
$23$ & $55$ & $0$ & $0$ & $\scalebox{.8}{5/24}$ \\
$24$ & $267$ & $2$ & $1$ & $\scalebox{.8}{1/8}$ \\
$25$ & $558$ & $10$ & $1$ & $\scalebox{.8}{1/12}$ \\
$26$ & $1\,888$ & $20$ & $0$ & $\scalebox{.8}{1/12}$ \\
$27$ & $3\,024$ & $29$ & $1$ & $0$ \\
\end{tabular} 
\end{center}
}
\caption{{\small Hunting ${\rm X}_{28}^R$ with $R=7{\bf A}_1\,3{\bf A}_2\,{\bf A}_7$ and $V=6{\bf A}_1\,2{\bf A}_2\,{\bf A}_7$.}}
\label{tab:7A13A2A7dim28}
\end{table}


\subsection{The root system $4{\bf A}_1\,2{\bf A}_2 \,2{\bf A}_3\,{\bf D}_4$ in dimension $28$}

${}^{}$\indent As another example with a little more lattices, as well as a component of type ${\bf D}$, consider the case of the root system $R$ above. 
The reduced mass of ${\rm X}_{28}^R$ is $1033/16$. We use the visible root system $V=3{\bf A}_1\,2{\bf A}_2 \,2{\bf A}_3\,{\bf D}_4$, that is $e=4$ and type $4^3\,3^2\,2^3\,1^4$ (hence $d$ is even $\geq 24$). The \texttt{BNE} algorithm terminates at $d=36$. It returns the $156$ lattices of ${\rm X}_{28}^R$, with reduced mass {\small $1/2$} ($112$ times), 
{\small $1/4$} ($26$ times), {\small $1/8$} ($15$ times) and {\small $1/16$} ($3$ times), after about $24$ h of CPU time. \ps

\tabcolsep=7pt
\begin{table}[H]
{\scriptsize 
\renewcommand{\arraystretch}{1.2} \medskip
\begin{center}
\begin{tabular}{c|c|c|c|c}
$d$ & $\sharp \texttt{iso}$ & $\sharp\texttt{found}$ & $\sharp \texttt{new\textunderscore lat}$ & $\texttt{rem\textunderscore red\textunderscore mass}$  \\
\hline
$24$ & $295$ & $78$ & $36$ & $\scalebox{.8}{773/16}$ \\
$26$ & $2\,082$ & $341$ & $66$ & $\scalebox{.8}{293/16}$ \\
$28$ & $12\,217$ & $1\,623$ & $44$ & $\scalebox{.8}{13/8}$ \\
$30$ & $55\,083$ & $4\,980$ & $6$ & $\scalebox{.8}{7/16}$ \\
$32$ & $154\,458$ & $10\,992$ & $4$ & $0$ \\
\end{tabular} 
\end{center}
}
\caption{{\small Hunting ${\rm X}_{28}^R$ with $R=4{\bf A}_1\,2{\bf A}_2 \,2{\bf A}_3\,{\bf D}_4$ and $V=3{\bf A}_1\,2{\bf A}_2 \,2{\bf A}_3\,{\bf D}_4$.}}
\label{tab:4A12A22A3D4dim28}
\end{table}


\subsection{The root system $16{\bf A}_1\, {\bf E}_6$ in dimension $28$}

${}^{}$ \indent 
For this root system $R$, the reduced mass of ${\rm X}_{28}^R$ is {\small $1/23040$}, hence there may well be a unique class in ${\rm X}_{28}^R$. 
The difficulty here is that no possible visible root system is ``very close'' to $R$, both because the presence of a component of type ${\bf E}$ and of too many ${\bf A}_1$. The best idea here is to use the fact that ${\bf A}_5 {\bf A}_1$ is a $2$-kernel of ${\bf E}_6$ (see Example 5.22 in \cite{cheuni}), and to run \texttt{BNE} for the visible root system $V \simeq 10{\bf A}_1\,{\bf A}_5\,{\bf D}_2 \simeq 12{\bf A}_1\,{\bf A}_5$, hence $e=2$ and type $6\,2^{11}$ (this forces $d$ even and $\geq 24$). For each of the even integers $24 \leq d \leq 46$,  it turns out that the values of $\sharp\texttt{iso}$ given by \texttt{BNE} are 
$$0, 0, 4, 5, 42, 93, 344, 516, 1440, 2064, 5792, 7673,$$ but that in all cases we have $\sharp\texttt{found}=0$. Nevertheless, for $d=48$ we find 
$\sharp\texttt{iso}=17098$ and finally $\sharp\texttt{found}=1$ ! The found (and unique) lattice in ${\rm X}_{28}^R$ has reduced mass {\small $1/23040$} and is ${\rm N}_{48}(x; 0)$ with
{\rm
$$x=(\scalebox{.7}{1, 1, 1, 1, 1, 1, 3, 3, 6, 6, 8, 8, 10, 10, 12, 12, 14, 14, 16, 16, 18, 18, 20, 20, 22, 22, 24, 24}) \in \Z^{28}.$$
}
The whole computation only takes about $13$ minutes of CPU time (essentially because we computed a single ${\rm BV}$ invariant). 
The explanation for all those zeros values of $\sharp\texttt{found}$ above is that this choice of visible root system, 
although essentially the best we can take, will most likely find lattices with root system of the form $n{\bf A}_1 {\bf A}_5$ with $n\geq 13$. 
Indeed, there are many lattices in ${\rm X}_{28}$ with root systems $13{\bf A}_1 {\bf A}_5$  (including some with the large reduced mass $1/4$) and $17{\bf A}_1{\bf A}_5$. 

\begin{remark} 
{\it
The root system ${\bf D}_5$ is also a $2$-kernel of ${\bf E}_6$, and the lattice $L$ above can also be found using the visible root system $11 {\bf A}_1 {\bf D}_5$ for $d=50$ (within about 55 minutes). In the latter case we have 
{\rm $\sharp\texttt{iso}=31790$, $\sharp\texttt{found}=7$,} and obtained $L \simeq {\rm N}_{50}(y)$ with $y \in \Z^{28}$ defined as
{\rm 
$$y=(\scalebox{.7}{1, 1, 5, 5, 6, 6, 7, 7, 8, 8, 12, 12, 14, 14, 15, 16, 16, 21, 21, 22, 22, 23, 23, 25, 25, 25, 25, 25}).$$
}
All other choices of visible root systems seem to require more computations.
}
\end{remark}

\subsection{The root system $8{\bf A}_1\,2{\bf A}_2$ in dimension $28$} 
\label{subsect:8A12A2}

${}^{}$ \indent For this root system $R$ the reduced mass of ${\rm X}_{28}^R$ is {\small $9694663/2880$} $\simeq 3366.2$, 
which already implies  $|{\rm X}_{28}^R|\geq 6733$. Actually, as we shall see, we have $|{\rm X}_{28}^R|=7603$, and this is actually the maximum of the $|{\rm X}_{28}^S|$ for all root systems $S$. This is one reason why we chosed this example, which is much harder than the previous ones. \ps
	We know from Sect. 5.9 \cite{cheuni} (and especially Example 5.14) that the pair $(R,V)$ is safe for $V=7{\bf A}_1\,2{\bf A}_2$, so we run $\texttt{BNE}$ for this visible rot system, that is $e\leq 1$ and type 
$3^2\,2^7\,1^8$	 ($s=17$). For the subsequent purpose of computing $2$-neighbors, it is convenient to restrict our search to odd $d$. Here is what $\texttt{BNE}$
	finds for all odd $d$ from $2s+1=35$ to $43$ (we add in the table a lower bound of the number $\sharp \texttt{rem\textunderscore lat}$ of remaining lattices in ${\rm X}_{28}^R$ after Step 7 in the last column):
\tabcolsep=7pt
\begin{table}[H]
{\scriptsize 
\renewcommand{\arraystretch}{1.2} \medskip
\begin{center}
\begin{tabular}{c|c|c|c|c|c}
$d$ & $\sharp \texttt{iso}$ & $\sharp\texttt{found}$ & $\sharp \texttt{new\textunderscore lat}$ & $\texttt{rem\textunderscore red\textunderscore mass}$ & $\sharp \texttt{rem\textunderscore lat}\geq$ \\
\hline
$35$ & $2\,039$ & $710$ & $518$ & $\scalebox{.8}{8987263/2880}$ & $6\,243$ \\
$37$ & $25\,009$ & $9\,587$ & $3\,536$ & $\scalebox{.8}{4199593/2880}$ & $2\,918$ \\
$39$ & $293\,217$ & $99\,546$ & $2\,958$ & $\scalebox{.8}{497041/2880}$ & $347$ \\
$41$ & $1\,280\,597$ & $367\,529$ & $451$ & $\scalebox{.8}{21349/720}$ & $61$ \\
$43$ & $5\,801\,141$ & $1\,398\,150$ & $43$ & $\scalebox{.8}{899/36}$ & $51$ \\
\end{tabular} 
\end{center}
}
\caption{
{\small Hunting ${\rm X}_{28}^R$ with $R=8{\bf A}_1\,2{\bf A}_2$ and  $V=7{\bf A}_1\,2{\bf A}_2$.}
}
\label{tab:8A12A2dim28}
\end{table}
This computation is already quite lengthy: it took about $1400$ h of CPU time.
It allowed to find $7506$ different classes in ${\rm X}_{28}^R$ so far, with remaining reduced mass $899/36$. 
It would be natural to try to go on in order to find the remaining lattices, namely by exploring $d=45$, $47$ and so on... 
This would unfortunately require much more CPU time. Instead, we may proceed as follows.\ps\ps

{\bf (a) Exceptional lattices in ${\rm X}_{28}^R$.} 
An inspection of the $7506$ found lattices is that none of them is exceptional. 
We refer to Lemma~\ref{lem:mincharvec} for an explanation of an analogous phenomenon for ${\rm X}_{29}^\emptyset$. 
Our aim now is to determine the (classes) of exceptional unimodular lattices in ${\rm X}_{28}^R$.\ps

As we are in dimension $28 \equiv 4 \bmod 8$, exceptional lattices $L$ in $\mathcal{L}_{28}$ with ${\rm r}_1(L)=0$
may be described using their singular companion $L':={\rm sing}(L)$, by Proposition~\ref{prop:singcomp}.
Write $L' \simeq {\rm I}_r \perp U$ with $r\geq 1$ and $U \in \mathcal{L}_{28-r}$ satisfying ${\rm r}_1(U)=0$. 
In particular, ${\bf D}_r \simeq {\rm R}_2({\rm I}_r)$ is an irreducible component of ${\rm R}_2(L')={\rm R}_2(L)$ (recall ${\bf D}_1=\emptyset$ and ${\bf D}_2\,=\,2{\bf A}_1$). Going back to the root system $R \simeq 8{\bf A}_1\,2{\bf A}_2$, 
the only possibilities are thus either $r=1$ and ${\rm R}_2(U) \simeq R$, or $r=2$ and ${\rm R}_2(U) \simeq 6{\bf A}_1\,2{\bf A}_2$. But an inspection of the classification in \cite{cheuni} of ${\rm X}_{26}$ and ${\rm X}_{27}$ shows that there are exactly two possibilities for $U$ for $r=2$, and $89$ for $r=1$ (all being non exceptional). 
Proposition~\ref{prop:singcomp} thus shows:

\begin{lemma} 
\label{lem:exc28R}
There are $89+2=91$ classes of exceptional lattices in ${\rm X}_{28}^R$. 
The sum of their reduced mass is {\small $595/24$} $+$ {\small $13/288$} $=$ {\small $7153/288$}.
\end{lemma}

Better, in the analysis above, if we have $U \,\simeq \,{\rm N}_d(x)$ with $d$ odd and $x \in \Z^{28-r}$, then we know from Lemma 11.1 and Proposition 11.6 in \cite{cheuni} that the corresponding $L$ satisfies $L \simeq {\rm N}_{2d}(y;0)$ where $y \in \Z^{28}$ has odd coordinates and satisfies $y_i \equiv x_i \bmod d$  for $i < 28-r$, and $y_i =d$ for $i\geq 28-r$. Neighbor forms for the $2+89$ different $U$ are easily found with $d$ odd using $\texttt{BNE}$ in dimensions $26$ and $27$. For instance, the $2$ classes $U$ in ${\rm X}_{26}$ with root system 
$6{\bf A}_1\, 2{\bf A}_2$ are those of the ${\rm N}_{35}(x)$ with $x \in \Z^{26}$ given by
$$x=\left\{ 
\begin{array}{l} 
\scalebox{.7}{(1, 1, 1, 2, 3, 4, 5, 6, 7, 8, 9, 10, 11, 11, 12, 12, 13, 13, 14, 14, 15, 15, 16, 17, 17, 17)},
\\ 
\scalebox{.7}{(1, 1, 1, 2, 3, 4, 5, 6, 7, 8, 9, 9, 10, 11, 11, 12, 12, 13, 13, 14, 14, 15, 16, 16, 16, 17)},
\end{array}\right.
$$
and their respective reduced mass is {\small $1/288$} and {\small $1/24$}. In the end, we do obtain neighbor forms
for all $91$ exceptional lattices of Lemma~\ref{lem:exc28R}. \ps
\ps

{\bf (b) Strict $2$ neighbors of rare lattices.} After incorporating the $91$ exceptional lattices above in ${\rm X}_{28}^R$, 
the remaining reduced mass is {\small $13/96$}. The remaining lattices presumably all have a small reduced mass. 
In order to find them, we apply some ideas from the theory of {\it visible isometries} explained in \cite[\S 7]{cheuni} (see especially \S 7.9):  \ps

(b1) The first idea is that if some $L \in \mathcal{L}_n$ has a large 
isometry group, we expect that many $2$-isotropic lines in 
$L/2L$ will be stable by a non trivial isometry of $L$, hence producing a 
$2$-neighbor having that isometry, and biasing the search. \ps

(b2) A second idea is that to in order to avoid computing all $2$-neighbors of $L$ (or ${\rm O}(L)$-orbits of them), it seems more promising to focus on those having the same visible root system as $L$ viewed as neighbors of ${\rm I}_n$; we call them the {\it strict} $2$-neighbors of $L$.
Concretely, if we have $L={\rm N}_d(x)$ with $d$ odd, $x \in \Z^n$ and $(d,x)$ normalized, then the strict $2$-neighbors of $L$ are the ${\rm N}_{2d}(y;\epsilon)$ with $y \in \Z^{n}$ satisfying $y \equiv x \bmod d$, as well as $y_1=1$ and, for all $1\leq i,j \leq n$ with $x_i = x_j$, the congruence $y_i \equiv y_j \bmod 2$. \ps

Among the $494$ lattices found above for $d=41$ and $d=43$, there are $12$ lattices with reduced mass $\leq 1/32$. 
For each of these $d$-neighbors $L$, given as $L={\rm N}_d(x)$, we compute all the {\it strict} $2$-neighbors $L'$ of $L$, 
with $L'={\rm N}_{2d}(y,\epsilon)$ as in (b2) above. 
As we are in dimension $n=28$ and $L$ has a visible root system $7{\bf A}_1\,2{\bf A}_2$, there are thus only $2^{28-2-7-4-1}=65\,536$ choices for $y$, hence presumably $32\,768$ isotropic ones, a quite manageable quantity. The hope is to find the remaining lattices among those $2\times 12\times 32\,768=786\,432$ different $2$-neighbors.  It works! Indeed, in about $20$ h of CPU time we did find this way the $6$ remaining elements of ${\rm X}_{28}^R$. They have a reduced mass $1/96$, $1/64$ twice and $1/32$ three times, and are respectively given by the following values of $(d,y,\epsilon)$:

\tabcolsep=2.5pt
\begin{table}[H]
{\scriptsize 
\renewcommand{\arraystretch}{1.2} \medskip
\begin{center}
\begin{tabular}{c|c|c|c }
$2d$ & $y$ & $\epsilon$ & $\texttt{red\textunderscore mass}$ \\
\hline
$82$ & $\scalebox{.7}{(1, 1, 1, 2, 3, 3, 4, 4, 36, 7, 33, 33, 31, 31, 11, 12, 13, 13, 14, 14, 14, 15, 25, 23, 22, 22, 20, 20)}$ & $0$ & $\scalebox{.8}{1/32}$ \\
$82$ & $\scalebox{.7}{(1, 1, 1, 2, 3, 3, 37, 37, 36, 7, 8, 8, 10, 10, 30, 29, 13, 13, 14, 14, 14, 15, 16, 18, 19, 19, 21, 21)}$ & $1$ & $\scalebox{.8}{1/96}$ \\
$86$ & $\scalebox{.7}{(1, 1, 1, 2, 2, 4, 4, 6, 6, 36, 9, 9, 33, 32, 32, 12, 30, 29, 29, 29, 28, 16, 26, 19, 19, 23, 22, 22)}$ & $0$ & $\scalebox{.8}{1/32}$ \\ 
$86$ & $\scalebox{.7}{(1, 1, 1, 2, 2, 4, 4, 6, 6, 36, 9, 9, 10, 32, 32, 12, 13, 29, 29, 29, 28, 16, 17, 19, 19, 20, 22, 22)}$ & $1$ & $\scalebox{.8}{1/32}$ \\
$86$ & $\scalebox{.7}{(1, 1, 1, 3, 3, 5, 37, 37, 36, 35, 9, 9, 11, 31, 31, 13, 14, 14, 14, 15, 16, 16, 26, 18, 18, 20, 22, 22)}$ & $1$ & $\scalebox{.8}{1/64}$ \\
$86$ & $\scalebox{.7}{(1, 1, 1, 39, 5, 5, 37, 7, 7, 35, 35, 10, 32, 32, 12, 13, 14, 14, 14, 28, 28, 27, 17, 17, 25, 24, 20, 20)}$ & $0$ & $\scalebox{.8}{1/64}$ 
\end{tabular}
\end{center}
}
\caption{
{\small The last ${\rm N}_{2d}(y;\epsilon)$ found in ${\rm X}_{28}^R$ for $R=8{\bf A}_1\,2{\bf A}_2$.}
}
\label{tab:8A12A2dim28bis}
\end{table}


\subsection{The empty root system in dimension $29$} 
\label{subsect:exdim29}

${}^{}$ \indent We finally discuss the determination of ${\rm X}_{29}^\emptyset$. 
For such lattices, the mass and the reduced mass coincide so we usually omit the term ``reduced''.
We know that the mass of ${\rm X}_{29}^\emptyset$
is\, {\small $49612728929/11136000$} $\simeq 4455.2$, so we have at least $8911$ isometry classes. 
The only possibility here is to take $V=\emptyset$, and we choose $0 \leq e \leq 1$ and all $n_i$ equal to $1$. 
Here is what $\texttt{BNE}$ finds for all {\it odd} $d$ from $2s+1=59$ to $83$, after about $850$ h of CPU time:

\tabcolsep=7pt
\begin{table}[H]
{\scriptsize 
\renewcommand{\arraystretch}{1.2} \medskip
\begin{center}
\begin{tabular}{c|c|c|c|c|c}
$d$ & $\sharp \texttt{iso}$ & $\sharp\texttt{found}$ & $\sharp \texttt{new\textunderscore lat}$ & $\texttt{rem\textunderscore red\textunderscore mass}$ & $\sharp \texttt{rem\textunderscore lat}\geq$ \\
\hline
$59$ & $1$ & $1$ & $1$ & $\scalebox{.8}{1710782101/384000}$ & $8\,912$ \\
$61 \,\&\, 63$ & $0$ & $0$ & $0$ & $\scalebox{.8}{1710782101/384000}$ & $8\,912$ \\
$65$ & $4$ & $4$ & $4$ & $\scalebox{.8}{570063367/128000}$ & $8\,909$ \\
$67$ & $19$ & $19$ & $19$ & $\scalebox{.8}{568975367/128000}$ & $8\,892$ \\
$69$ & $149$ & $138$ & $107$ & $\scalebox{.8}{562979367/128000}$ & $8\,798$ \\
$71$ & $654$ & $598$ & $527$ & $\scalebox{.8}{530939367/128000}$ & $8\,297$ \\
$73$ & $3\,173$ & $2\,771$ & $1\,836$ & $\scalebox{.8}{1254618101/384000}$ & $6\,536$ \\
$75$ & $24\,641$ & $20\,300$ & $3\,971$ & $\scalebox{.8}{185034167/128000}$ & $2\,893$ \\
$77$ & $70\,121$ & $55\,094$ & $2\,774$ & $\scalebox{.8}{276881503/1152000}$ & $482$ \\
$79$ & $206\,343$ & $153\,700$ & $605$ & $\scalebox{.8}{40169503/1152000}$ & $71$ \\
$81$ & $1\,029\,214$ & $725\,560$ & $96$ & $\scalebox{.8}{26594503/1152000}$ & $48$ \\
$83$ & $2\,321\,088$ & $1\,548\,714$ & $24$ & $\scalebox{.8}{25345453/1152000}$ & $46$ \\
\end{tabular} 
\end{center}
}
\caption{
{\small Hunting ${\rm X}_{29}^\emptyset$ with $V=\emptyset$ and odd $59 \leq d\leq 83$.}
}
\label{tab:X29emptyset}
\end{table}
So far we have found $9964$ elements in ${\rm X}_{29}^\emptyset$.  
For instance, the first of theses lattices, found for $d=59$, is the lattice ${\rm N}_{59}(1,2,3,\dots,29)$ 
which incidentally belongs to the family studied in~\cite[\S 8]{cheuni}. 
Also, an inspection of our list shows that we only found a single exceptional lattice so far, for $d=83$.
This can be partially explained by the following lemma.

\begin{lemma} 
\label{lem:mincharvec}
Assume $L \in \mathcal{L}_n$ is a $p$-neighbor of ${\rm I}_n$ with $p$ prime and empty visible root 
system {\rm (}this forces $p\geq 2n+1${\rm )}. Then any characteristic vector $\xi$ of $L$ with $\xi.\xi<n$
satisfies $\xi.\xi \geq \frac{4n^3-n}{3p^2}$. \end{lemma}

\begin{pf} 
By definition, we have $L={\rm M}_p(x) + \Z \frac{x'}{p}$ for some $p$-isotropic $x \in \Z^n$
and some $x' \in \Z^n$ satisfying $x' \equiv x \bmod p$.  
As $p$ is odd, the vector $p \xi \in \Z^n$ is a characteristic vector of ${\rm I}_n$, hence has odd coordinates.
Write $p\xi = kx +p m$ with $0 \leq k <p$ and $m \in \Z^n$. 
We have $k>0$, otherwise $\xi \in \Z^n$ and $\xi.\xi \geq n$.
So the coordinates of $p \xi$ are distinct mod $p$, 
since so are those of $x$ by assumption. This proves 
$p^2 \xi . \xi \geq \sum_{i=1}^n (2i-1)^2 = (4n^3-n)/3$. \end{pf}

Note that for $n=29$ and $\xi.\xi=5$, this forces $p\geq 83$, in accordance with what we found.
Our first aim now is to seek for exceptional lattices in ${\rm X}_{29}^\emptyset$.
We cannot argue as in \S~\ref{subsect:8A12A2} (a) since $29 \not \equiv 4 \bmod 8$. 
Instead, we use the variant of \texttt{BNE} discussed in~\cite[\S 9.13--9.16]{cheuni}.
The basic idea is to look for $d$-neighbors $N:={\rm N}_d(x;\epsilon)$ of ${\rm I}_{29}$ with empty root system and such that 
the norm $5$ vector $(0,\dots,0,1,1,1,1,1) \in {\rm I}_{29}$ is a (visible!) characteristic vector of $N$.
As explained {\it loc. cit.}, the trick is just to modify Step 2 of \texttt{BNE} and enumerate only $d$-isotropic $x \in \Z^{29}$, with $d$ even, 
all of whose coordinates are odd, except the last $5$ ones which are even and with sum $\equiv 0 \bmod d$. This forces $d\geq 94$,
and we obtain after about $60$ h of CPU time:

\tabcolsep=7pt
\begin{table}[H]
{\scriptsize 
\renewcommand{\arraystretch}{1.2} \medskip
\begin{center}
\begin{tabular}{c|c|c|c|c|c}
$d$ & $\sharp \texttt{iso}$ & $\sharp\texttt{found}$ & $\sharp \texttt{new\textunderscore lat}$ & $\texttt{rem\textunderscore mass}$ & $\sharp \texttt{rem\textunderscore lat}$  \\
\hline
$94$ & $20$ & $20$ & $7$  & $\scalebox{.8}{22177453/1152000}$ & $40$\\
$96$ & $46$ & $37$ & $9$ & $\scalebox{.8}{2697341/144000}$ & $39$ \\
$98$ & $82$ & $74$ & $26$ & $\scalebox{.8}{1488341/144000}$ & $22$ \\
$100$ & $150$ & $122$ & $24$ & $\scalebox{.8}{729281/144000}$ & $12$   \\
$102$ & $900$ & $664$ & $26$  & $\scalebox{.8}{123131/144000}$ & $3$ \\
$104$ & $687$ & $466$ & $6$ & $\scalebox{.8}{81431/144000}$ & $3$ \\
$106$ & $4\,940$ & $3\,131$ & $4$  & $\scalebox{.8}{40631/144000}$ & $2$ \\
$108$ & $7\,833$ & $4\,627$ & $1$  & $\scalebox{.8}{37031/144000}$ & $2$\\
$110$ & $55\,116$ & $29\,680$ & $0$  & $\scalebox{.8}{37031/144000}$ & $2$ \\
$112$ & $47\,310$ & $23\,889$ & $0$  & $\scalebox{.8}{37031/144000}$ & $2$ \\
$114$ & $377\,410$ & $176\,966$ & $1$  & $\scalebox{.8}{1481/5760}$ & $2$ \\
\end{tabular} 
\end{center}
}
\caption{
{\small Hunting exceptional lattices in ${\rm X}_{29}^\emptyset$.}
}
\label{tab:X29emptysetex}
\end{table}
The last lattice, with mass {\small $1/24000$}, came very late and is much harder to find than the others: this is ${\rm N}_{114}(x;0)$ with $$x=\scalebox{.7}{(1, 2, 3, 5, 7, 9, 11, 15, 16, 17, 21, 22, 23, 27, 29, 31, 33, 35, 36, 37, 38, 39, 41, 45, 49, 51, 53, 55, 57)}.$$
Up to this point, we have found $10068$ classes in ${\rm X}_{29}^\emptyset$, with remaining mass {\small $1481/5760$}.
There are several methods that we can use to find the last ones. First, we run \texttt{BNE}
for the empty visible root system, $0\leq e \leq 1$, and {\it even} $58 \leq d \leq 82$, which we have not done yet, and only finds $15$ more lattices in about $1050$ more hours:

\tabcolsep=7pt
\begin{table}[H]
{\scriptsize 
\renewcommand{\arraystretch}{1.2} \medskip
\begin{center}
\begin{tabular}{c|c|c|c|c|c}
$d$ & $\sharp \texttt{iso}$ & $\sharp\texttt{found}$ & $\sharp \texttt{new\textunderscore lat}$ & $\texttt{rem\textunderscore red\textunderscore mass}$ & $\sharp \texttt{rem\textunderscore lat}\geq$ \\
\hline
$58 - 70$ & $1\,207$ & $1\,112$ & $0$ & $\scalebox{.8}{1481/5760}$ & $2$ \\
$72$ & $4\,445$ & $3\,847$ & $2$ & $\scalebox{.8}{1091/5760}$ & $2$ \\
$74$ & $16\,304$ & $13\,423$ & $0$ & $\scalebox{.8}{1091/5760}$ & $2$ \\
$76$ & $65\,591$ & $51\,666$ & $3$ & $\scalebox{.8}{589/3840}$ & $2$ \\
$78$ & $390\,922$ & $290\,029$ & $2$ & $\scalebox{.8}{499/3840}$ & $2$ \\
$80$ & $1\,065\,081$ & $752\,393$ & $7$ & $\scalebox{.8}{941/15360}$ & $2$ \\
$82$ & $2\,969\,999$ & $1\,985\,285$ & $1$ & $\scalebox{.8}{301/15360}$ & $2$ \\
\end{tabular} 
\end{center}
}
\caption{
{\small Hunting ${\rm X}_{29}^\emptyset$ with $V=\emptyset$ and even $58 \leq d\leq 82$.}
}
\label{tab:X29emptyseteven}
\end{table}
At this point, we have found $10\,083$ lattices, and the remaining mass is only {\small $301/15360$}.
It would be presumably quite lengthy to seek for the remaining lattices by simply pursuing \texttt{BNE}. 
A more direct way to find them is to use the theory of visible isometries as in \S\ref{subsect:8A12A2} (b),
and study $2$-neighbors of the found lattices $L={\rm N}_d(x)$ with large isometry groups and $d$ odd. 
Note that as the root system is empty here, each $2$ neighbor of such an $L$ is strict, 
so we have no meaningful way to reduce the search as we did in \S\ref{subsect:8A12A2}.
In such a situation with ${\rm O}(L)$ big, we should 
also gain much in principle in computing first the ${\rm O}(L)$-orbits of 
$2$-isotropic lines in $L/2L$, since there is a huge number of such lines in dimension $29$. 
In practice, this is not really necessary, and 
we prefer to only compute the neighbors associated to a large number 
of $2$-isotropic lines of a given $L$, say here $2\,000\,000$, and then try the 
next $L$ if it fails.\ps

This works pretty well! Indeed, the lattice $L={\rm N}_d(x)$ in ${\rm X}_{29}^\emptyset$, with $d$ odd and largest
isometry group, appears for $d=83$ and satisfies $|{\rm O}(L)|=1536$. Using this $L$ as explained above 
we do find $7$ new lattices, with remaining mass $167/25920$. The next $3$ lattices with largest isometry groups did not provide new lattices, but the $2$ after, namely a $75$-neighbor with mass $1/160$, and a $81$-neighbor with mass $1/128$, do give rise (each) to a new lattice, with respective masses {\small $1/160$} and {\small $1/5184$}, and concludes the proof! 
Those $7+1+1=9$ lattices ${\rm N}_{2d}(y;\epsilon)$ are given in the following table. The full computation here 
took about $480$ h of CPU time.

\tabcolsep=2.5pt
\begin{table}[H]
{\scriptsize 
\renewcommand{\arraystretch}{1.2} \medskip
\begin{center}
\begin{tabular}{c|c|c|c}
$2d$ & $y$ & $\epsilon$ & $\texttt{mass}$ \\
\hline
$166$ & $\scalebox{.7}{(1, 5, 6, 11, 12, 13, 15, 16, 65, 19, 20, 21, 23, 58, 57, 27, 28, 29, 53, 31, 51, 33, 49, 35, 36, 45, 44, 43, 42)}$ & $1$ & $\scalebox{.8}{1/864}$ \\
$166$ & $\scalebox{.7}{(1, 5, 6, 11, 12, 13, 15, 67, 18, 64, 20, 62, 23, 58, 57, 27, 55, 54, 30, 31, 32, 33, 49, 35, 36, 38, 39, 43, 42)}$ & $0$ & $\scalebox{.8}{1/1536}$ \\
$166$ & $\scalebox{.7}{(1, 5, 6, 11, 12, 13, 15, 16, 18, 19, 20, 62, 23, 58, 26, 27, 28, 29, 30, 52, 32, 50, 49, 48, 47, 38, 39, 40, 42)}$ & $1$ & $\scalebox{.8}{1/6144}$\\
$166$ & $\scalebox{.7}{(1, 5, 6, 11, 12, 13, 15, 16, 18, 64, 20, 21, 60, 58, 26, 56, 55, 54, 30, 31, 51, 50, 49, 35, 47, 38, 44, 40, 42)}$ & $0$ & $\scalebox{.8}{1/96}$ \\
$166$ & $\scalebox{.7}{(1, 5, 6, 11, 12, 13, 15, 16, 65, 64, 20, 21, 23, 25, 57, 27, 55, 54, 53, 52, 32, 50, 49, 48, 36, 45, 39, 40, 42)}$ & $0$ & $\scalebox{.8}{1/2592}$\\
$166$ & $\scalebox{.7}{(1, 5, 6, 11, 12, 13, 15, 16, 18, 19, 63, 62, 23, 25, 26, 56, 28, 29, 53, 31, 51, 50, 34, 48, 47, 45, 39, 40, 42)}$ & $1$ & $\scalebox{.8}{1/3072}$ \\
$166$ & $\scalebox{.7}{(1, 5, 6, 11, 12, 13, 68, 16, 65, 64, 20, 21, 60, 58, 57, 56, 28, 29, 53, 31, 51, 33, 34, 35, 47, 45, 39, 40, 42)}$ & $0$ & $\scalebox{.8}{1/18432}$ \\
$150$ & $\scalebox{.7}{(1, 4, 5, 6, 7, 9, 65, 11, 12, 62, 61, 15, 59, 18, 56, 55, 22, 24, 25, 49, 48, 28, 46, 44, 43, 41, 35, 39, 38)}$ & $1$ & $\scalebox{.8}{1/160}$ \\
$162$ & $\scalebox{.7}{(1, 5, 8, 9, 11, 12, 13, 14, 15, 65, 17, 18, 19, 20, 60, 22, 57, 55, 27, 52, 30, 49, 33, 47, 46, 36, 38, 42, 41)}$ & $1$ & $\scalebox{.8}{1/5184}$
\end{tabular} 
\end{center}
}
\caption{
{\small The last nine lattices ${\rm N}_{2d}(y;\epsilon)$ found in ${\rm X}_{29}^\emptyset$}
}
\label{tab:X29empty2nei}
\end{table}


The method above is essentially the one we used when we first computed ${\rm X}_{29}^\emptyset$. 
Meanwhile, the second author found a more direct way to find $d$-neighbors of ${\rm I}_n$ with 
prescribed (and ``visible'') isometries. This method is described in \S 7.5 and \S 7.7 of \cite{cheuni}. 
The basic idea is to fix $\sigma \in {\rm O}({\rm I}_n)$ and to study to $d$-neighbors $N$ of ${\rm I}_n$ with 
$\sigma \in {\rm O}(N)$ and having a given visible root system $V$, which translates into some conditions on the isotropic lines 
that we enumerate. This is especially suited to empty (or small) $V$.
An example of application of these ideas to the determination of ${\rm X}_{28}^\emptyset$ is detailed in \S 7.6 {\it loc. cit.}.
The situation here is a bit similar, and goes as follows.\ps
 We go back right before the $2$-neighbor argument above. At this step the remaining mass is {\small $301/15360$}.
 We have $15\,360 \,=\, 2^{10}\,3\,5$, so we know that we still need to find lattices in ${\rm X}_{29}^\emptyset$
 having isometries of prime order $q=2, 3$ and $5$. The characteristic polynomial of such an isometry is 
 $\phi_q^k \,\phi_1^l$, with $\phi_m$ the $m$-th cyclotomic polynomial, and $k(q-1)+l=29$. 
 For fixed $q$ and $k$, we choose an auxiliary odd prime $p \equiv 1 \bmod q$ and consider associated $d:=pd'$ isotropic 
 lines describes in {\it loc. cit.} \S 7.5 for all odd integers $d'=1,3,\dots$. 
 These lines have the properties that the associated neighbors $N$ have an empty root system and a (visible) isometry 
 with characteristic polynomial $\phi_q^k \,\phi_1^l$. Better, they are tailored such that 
 the following order $q$ element 
$$\sigma_{q,k}\,\,=\,\,(1\,2\,\dots \,\,q)\,\,(q+1\,\,\,q+2\,\,\,\dots\,\,\,2q)\,\, \cdots \,\,(\,(k-1)q+1\,\,\,(k-1)q+2\,\,\,\dots \,\,\,kq),$$
a product of $k$\,\,disjoint cycles  in ${\rm S}_{29}$, lies in ${\rm O}(N)$. This requires $qk \leq 29$. 
By studying those lines we do also find the $9$ remaining lattices, under the form given by Table~\ref{tab:X29empty2neibis} below. This is much faster: it only took less than $3$ h to find the first $8$ lattices, and about $7$ h for the last one.
 

\hspace{-2cm}
\tabcolsep=2pt
\begin{table}[H]
{\scriptsize 
\renewcommand{\arraystretch}{1.2} \medskip
\begin{center}
\begin{tabular}{c|c|c|c|c|c}
$d$ & $x$ & $\texttt{mass}$ & $\texttt{char}$ & $p$ & $d'$ \\
\hline
$407$ & $\scalebox{.5}{(1, 334, 38, 75, 223, 78, 4, 115, 152, 300, 375, 301, 5, 42, 190, 266, 118, 340, 7, 303, 45, 378, 82, 119, 267, 231, 121, 11, 308)}$ & $\scalebox{.8}{1/160}$ & $\phi_5^5\phi_1^{9}$ & $11$ & $37$ \\
$315$ & $\scalebox{.5}{(1, 226, 46, 92, 2, 137, 183, 93, 228, 94, 139, 229, 50, 275, 95, 141, 51, 186, 8, 233, 53, 190, 100, 235, 56, 147, 238, 14, 105)}$ & $\scalebox{.8}{1/2592}$ & $\phi_3^8\phi_1^{13}$ & $7$ & $45$ \\
$315$ & $\scalebox{.5}{(1, 226, 46, 92, 2, 137, 3, 48, 138, 274, 184, 4, 50, 275, 95, 276, 6, 96, 8, 233, 53, 190, 100, 235, 56, 147, 238, 14, 105)}$ & $\scalebox{.8}{1/864}$ & $\phi_3^8\phi_1^{13}$ & $7$ & $45$ \\
$357$ & $\scalebox{.5}{(1, 205, 256, 309, 156, 207, 106, 310, 4, 311, 209, 5, 211, 58, 109, 213, 111, 264, 164, 62, 215, 64, 268, 319, 14, 168, 322, 119, 273)}$ & $\scalebox{.8}{1/5184}$ & $\phi_3^8\phi_1^{13}$ & $7$ & $51$ \\
$287$ & $\scalebox{.5}{(1, 247, 165, 248, 125, 166, 3, 167, 208, 127, 86, 4, 169, 128, 46, 253, 212, 130, 8, 254, 172, 50, 9, 214, 10, 174, 215, 175, 217)}$ & $\scalebox{.8}{1/3072}$ & $\phi_3^9\phi_1^{11}$ & $7$ & $41$ \\
$287$ & $\scalebox{.5}{(1, 247, 165, 248, 125, 166, 3, 167, 208, 45, 209, 250, 169, 128, 46, 253, 212, 130, 213, 90, 131, 50, 9, 214, 10, 174, 215, 175, 217)}$ & $\scalebox{.8}{1/96}$ & $\phi_3^9\phi_1^{11}$ & $7$ & $41$ \\
$301$ & $\scalebox{.5}{(1, 44, 130, 260, 2, 88, 218, 261, 46, 134, 177, 263, 178, 6, 264, 136, 265, 222, 267, 9, 95, 225, 268, 53, 183, 226, 11, 98, 56)}$ & $\scalebox{.8}{1/6144}$ & $\phi_3^9\phi_1^{11}$ & $7$ & $43$ \\
$329$ & $\scalebox{.5}{(1, 142, 95, 143, 237, 96, 239, 51, 4, 99, 240, 193, 288, 100, 53, 148, 289, 242, 197, 9, 291, 199, 293, 152, 295, 107, 60, 14, 203)}$ & $\scalebox{.8}{1/18432}$ & $\phi_3^9\phi_1^{11}$ & $7$ & $47$ \\
$483$ & $\scalebox{.5}{(1, 139, 416, 347, 3, 417, 73, 4, 352, 145, 423, 354, 79, 355, 289, 82, 359, 152, 429, 222, 16, 292, 17, 431, 294, 434, 21, 91, 161)}$ & $\scalebox{.8}{1/1536}$ & $\phi_2^{12}\phi_1^{17}$ & $7$ & $67$ \\
\end{tabular} 
\end{center}
}
\caption{
{\small Another form for the last nine lattices ${\rm N}_d(x)$ found in ${\rm X}_{29}^\emptyset$}
}
\label{tab:X29empty2neibis}
\end{table}

For instance, for the lattice $L={\rm N}_{329}(x)$ above with mass $1/18432$ we have 
$$\begin{array}{ccc} 
x \bmod 7 & = & \scalebox{.7}{(1, 2, 4, 3, 6, 5, 1, 2, 4, 1, 2, 4, 1, 2, 4, 1, 2, 4, 1, 2, 4, 3, 6, 5, 1, 2, 4, 0, 0)}, \\
x \bmod 47 & = & \scalebox{.7}{(1, 1, 1, 2, 2, 2, 4, 4, 4, 5, 5, 5, 6, 6, 6, 7, 7, 7, 9, 9, 9, 11, 11, 11, 13, 13, 13, 14, 15),}
\end{array}
$$
so $\sigma_{3,9}^{-1}$ acts on $(\Z/7)x  \subset (\Z/7)^{29}$ by multiplication by $2$, and fixes $x \bmod 47$.

\begin{remark} 
This second method also has drawbacks.
As an example, consider the last lattice $L$ in Table~\ref{tab:X29empty2neibis}. 
Although we have $|{\rm O}(L)|=1536 = 2^9\, \cdot 3$, 
it is impossible to find it as a neighbor of ${\rm I}_{29}$ with a visible isometry of order $3$. 
Indeed, we can check a posteriori that the order $3$ elements of ${\rm O}(L)$ all have the same characteristic polynomial $\phi_3^{11} \phi_1^7$, whereas no element of ${\rm O}({\rm I}_{29})$ has this property since $11 \cdot 3 > 29$. They are other constraints, which luckily are not prohibitive to conclude here: see~\S 7.7 {\it loc. cit.} for more about this. 
\end{remark}

\end{document}